\documentclass{article}

\usepackage[english]{babel}
\usepackage[style=numeric]{biblatex}
\usepackage[letterpaper,top=2cm,bottom=2cm,left=3cm,right=3cm,marginparwidth=1.75cm]{geometry}

\usepackage{amsmath, amssymb}
\usepackage{amsthm}
\usepackage{graphicx}
\usepackage{hyperref}
\hypersetup{colorlinks=true,allcolors=blue}
\usepackage{soul, color}
\usepackage{tabularx}
\usepackage{booktabs}
\usepackage{csquotes}
\usepackage{algorithm}
\usepackage{authblk}
\usepackage{algpseudocode}
\usepackage{multirow}
\usepackage{longtable}
\usepackage{csvsimple-l3}

\title{Solving Highly Constrained Multi-Objective Decision Problems\\
with Preference-Guided Metaheuristics:\\
\large The IMAP-IGS Framework}
\author[1]{Timp, L.J.K. (Lennard)}
\author[1]{Yorke-Smith, N. (Neil)}
\author[1]{Wolfert, A.R.M. (Rogier)}

\affil[1]{Delft University of Technology}

\date{\today}
\begin{document}

    \maketitle

\begin{abstract}
    Highly constrained multi-objective design and decision problems arise throughout operations research and remain difficult to solve, due to tight feasibility requirements and conflicting stakeholder preferences. Population-based evolutionary algorithms are the standard heuristic tool for this problem class, because they search large, non-convex, and combinatorial solution spaces without requiring gradients or convexity. Existing evolutionary approaches, however, typically treat constraint handling and multi-objective optimisation as separate mechanisms. Pareto-based methods provide a useful dominance ordering, but they are not inherently designed to identify a single group-preference optimum; they generate a Pareto front that requires a posteriori selection. Moreover, this post-processing step may fail to recover the true best-fit design--decision solution under strong preference conflict. Conversely, classical scalarisation methods collapse heterogeneous preference structures into a single objective, which can introduce structural distortions under multi-actor settings and obscure meaningful trade-offs.

    To address these limitations, this paper introduces the Integrative Maximisation of Aggregated Preferences (IMAP) method and its Inter-Generational Solver (IMAP-IGS), a metaheuristic framework that embeds IMAP scoring directly into the evolutionary fitness signal, thereby inducing a preference-consistent selection structure over the design--decision space. A constraint-violation preference function is introduced to steer the search toward feasibility within the same preference-based fitness structure, while final reported decision scores remain based only on substantive stakeholder preferences. Any population-based evolutionary algorithm that relies on scalar fitness evaluation can adopt IMAP-IGS by replacing its fitness function, while preserving all underlying evolutionary operators. A current-best re-insertion mechanism is introduced to stabilise intergenerational dynamics arising from z-normalised preference evaluation within the evolutionary loop.

    Two complementary instantiations are presented: IMAP-BRKGA, which embeds IMAP into the Biased Random-Key Genetic Algorithm (BRKGA) for combinatorial problems with hard feasibility structures, and IMAP-GA-II, a real-coded genetic algorithm that retains standard SBX crossover and polynomial mutation while removing Pareto-specific mechanisms such as non-dominated sorting and crowding-distance selection from NSGA-II.

    The frameworks are evaluated on three structurally distinct, highly constrained benchmark classes: the DAS-CMOP suite (nine problems, 16 difficulty levels, 30 seeds each), the MO-VRPTW Solomon instances (56 problems), and the Heterogeneous Vessel Allocation and Scheduling Problem (HVASP, ten synthetic instances modelled on real-world offshore operations, in two size classes). Under conflicting stakeholder preferences, IMAP-BRKGA consistently outperforms its weighted-sum BRKGA baseline in combinatorial settings, while IMAP-GA-II outperforms state-of-the-art Pareto-based CMOEA methods on DAS-CMOP. When preferences are aligned, these advantages diminish significantly, confirming that the observed gains stem from preference-conflict resolution rather than generic optimisation superiority.

    These results suggest that preference aggregation grounded in preference function modelling is particularly effective in the regime where real-world constrained design--decision problems are most difficult: settings with tightly coupled feasibility constraints and genuinely conflicting stakeholder preferences.
\end{abstract}


\section{Introduction}
\label{sec:introduction}

Highly constrained multi-objective design--decision problems are among the hardest practical problem classes in operations research. They combine tight feasibility requirements with competing performance dimensions and stakeholder preferences that may point in different directions. Vehicle routing under strict time windows, offshore vessel allocation, and multi-period resource scheduling all require the same kind of output: not a set of interesting alternatives, but a single implementable decision that best fits the group preference structure.

Existing evolutionary approaches usually separate three parts of this task. Constraint handling is treated through constraint domination, penalties, repair, or decoders. Multi-objective search is treated through Pareto dominance, decomposition, or scalarisation. Preference-based decision-making is often applied only after search, when the decision-maker selects one representative from the generated set. This separation is effective for exploratory analysis, but it is structurally indirect when the problem is already formulated as a design--decision problem with stakeholder preference functions and weights. Pareto-based methods provide a useful dominance ordering and can expose trade-offs, but they are not inherently designed to produce a unique group-preference optimum. Classical scalarisation produces a single search signal, but it can collapse heterogeneous preference structures into a single objective in a way that obscures stakeholder conflict.

The distinction is not whether Pareto search is useful. It is useful when the task is to order or explore alternatives. The distinction is whether the algorithmic search signal is aligned with the decision question. If the final question is ``which feasible solution has the highest aggregated group preference?'', then the search process should ideally use the same preference structure that defines the decision. A Pareto-based CMOEA searches for efficient trade-offs in objective space and then requires a selection rule. IMAP-IGS instead searches directly for the decision point induced by the supplied preference functions.

Preference Function Modelling (PFM) provides the measurement-theoretic foundation for this claim \cite{Barzilai2010}. In the ODESYS/IMAP line of work, decision validity requires that system performance dimensions are translated into a common preference domain, that system capability and stakeholder desirability are integrated, that individual stakeholder preferences can be associated and weighted, and that the method yields a unique best-fit decision \cite{Van_Heukelum2023Socio-technical,wolfert2026preferencebasedoptimisationgroupdecisionmaking}. The z-normalised weighted-centroid aggregator used by IMAP follows from this preference-aggregation argument: stakeholder preference scores are first expressed on affine preference scales and then aggregated as a weighted centroid in a common linear preference space \cite{wolfert2026uniqueaggregation}. These references establish the theoretical basis for preference aggregation; the present paper studies how that preference signal can be embedded directly into evolutionary search.

This paper introduces the Integrative Maximisation of Aggregated Preferences Inter-Generational Solver (IMAP-IGS), a host-agnostic metaheuristic framework for constrained multi-objective design--decision problems. The framework replaces the host algorithm's scalar fitness signal with IMAP scoring while preserving the host's representation, crossover, mutation, decoder, and other variation operators. This makes the framework applicable beyond one specific GA architecture. We demonstrate this in two structurally different instantiations. IMAP-BRKGA embeds IMAP into the Biased Random-Key Genetic Algorithm \cite{goncalves2013brkga}, whose decoder-based structure is suited to combinatorial problems with hard feasibility requirements. IMAP-GA-II is a real-coded GA instantiation that retains simulated binary crossover (SBX) and polynomial mutation while removing Pareto-specific mechanisms such as non-dominated sorting and crowding-distance selection from NSGA-II.

Two methodological additions make the embedding operational in highly constrained settings. First, a constraint-violation preference function is introduced to steer the search toward feasibility within the same preference-based fitness structure. It is a search-time pseudo-preference, not a substantive stakeholder objective, and it is excluded from final reported decision scores. Second, current-best re-insertion stabilises the intergenerational dynamics created by z-normalised IMAP evaluation. Because IMAP scores are computed relative to the current evaluation pool, the best currently known feasible solution must be repeatedly reintroduced into the live pool before scoring so that it remains comparable to the evolving population.

The paper's empirical claim is that IMAP-IGS is especially effective where feasibility constraints and preference conflict are coupled. This claim is tested across three structurally distinct benchmark classes under two stakeholder configurations. Configuration C creates conflicting preferences; Configuration A aligns preferences. If the benefit comes from resolving preference conflict, it should be strong in Configuration C and much smaller in Configuration A. A supplementary many-objective experiment (Section~\ref{subsec:results-dascmaop}) additionally verifies that the framework carries over unchanged as the number of objectives grows.

The contributions of this paper are:

\begin{enumerate}
    \item \textbf{IMAP-IGS framework:} a host-agnostic metaheuristic framework in which IMAP scoring replaces the scalar fitness signal of a population-based evolutionary algorithm while preserving the host algorithm's variation operators.
    \item \textbf{Constraint-violation preference function:} a search-time pseudo-preference function that steers the search toward feasibility inside the IMAP-IGS fitness structure, while final reported scores remain based only on substantive stakeholder preferences.
    \item \textbf{Current-best re-insertion:} a stabilisation mechanism for z-normalised intergenerational preference scoring that repeatedly reintroduces the best currently known feasible IMAP solution into the live evaluation pool.
    \item \textbf{Two algorithmic instantiations:} IMAP-BRKGA for combinatorial problems with hard feasibility structures and IMAP-GA-II for real-coded constrained benchmark problems.
    \item \textbf{Empirical validation under preference conflict:} evidence across the Difficulty Adjustable and Scalable CMOP (DAS-CMOP) benchmark suite \cite{Fan_2020_DASCMOP}, the Multi-Objective Vehicle Routing Problem with Time Windows (MO-VRPTW) Solomon benchmark \cite{solomon1987VRPTW}, and the Heterogeneous Vessel Allocation and Scheduling Problem (HVASP) showing that the advantage of IMAP-IGS is strongest when stakeholder preferences genuinely conflict and diminishes when preferences are aligned.
\end{enumerate}

The paper is structured as follows. Section~\ref{sec:literature_review} surveys the relevant literature and formalises the research gap. Section~\ref{sec:problem_formulation} defines the problem formulation and benchmark classes. Section~\ref{sec:methodology} presents the IMAP-IGS framework and both instantiations. Section~\ref{sec:experiments} describes the experimental setup. Section~\ref{sec:results} reports the results. Section~\ref{sec:discussion} discusses scope, limitations, and future work. Section~\ref{sec:conclusion} concludes.

\section{Literature Review}
\label{sec:literature_review}

\subsection{Highly Constrained Multi-Objective Decision Problems}
Highly constrained multi-objective decision problems arise wherever hard operational requirements coexist with competing performance objectives and genuine disagreement among stakeholders about what constitutes a good solution. In logistics, carriers must route fleets under hard time-window and capacity constraints while simultaneously minimising cost, travel time, and fleet size for customers with different service priorities. In offshore project scheduling, vessels must be assigned to activities under strict eligibility, availability, and non-overlap constraints \cite{Guo2023Integrated, Homsi2020Industrial} while minimising total mobilisation cost and project duration for an operator and a contractor whose objectives partially conflict \cite{Teuber2024}. In manufacturing, production schedules must satisfy machine availability, sequencing, and due-date constraints while balancing throughput, tardiness, and energy use across product lines with competing profitability profiles \cite{ghasemi2024}.

What distinguishes this class of problems from standard multi-objective optimisation is not the presence of constraints per se, constraint handling is a mature research area \cite{coello2002}, but the combination of three properties that existing algorithms rarely address jointly. First, the feasible region may be non-convex, disconnected, or extremely small relative to the decision space, making feasibility itself a hard subproblem. Second, the performance dimensions are genuinely conflicting across stakeholders, so no single efficient solution is universally preferred. Third, decision-makers require a single implementable solution, not a set of options for further deliberation. Together, these properties mean that the standard CMOEA pipeline, which converges to a Pareto-front approximation and then selects, is indirect for the problem considered here: it may be unreliable in the presence of tight feasibility constraints, and it defers the preference-sensitive part of the decision to a posteriori selection rather than using the preference model during search.

\subsection{Constrained Multi-Objective Evolutionary Algorithms}
The dominant paradigm in evolutionary multi-objective optimisation is Pareto-based selection. NSGA-II \cite{deb2002fast} uses fast non-dominated sorting and crowding distance to maintain a diverse population near the Pareto front. NSGA-III \cite{deb2014evolutionary} extends this to many-objective problems using reference vectors. The Constrained Two-Archive Evolutionary Algorithm (C-TAEA) \cite{li2019ctaea} maintains two co-evolving archives, one convergence-focused and one diversity-focused, to improve coverage of the Pareto front in constrained settings. The Multi-Objective Evolutionary Algorithm based on Decomposition with Constraint Dominance Principle (MOEA/D-CDP) \cite{Zhang2007MOEAD,Fan2016moea/d-cdp} decomposes the problem into scalar subproblems using a neighbourhood-based constraint-domination principle.

All of these methods share the same structural limitation for the design--decision problem class considered here: they produce a Pareto-front approximation rather than a preferred decision. The final step, selecting one solution from the front, is left to the decision-maker or delegated to an a posteriori mechanism such as weighted utility scoring or visual inspection. Pareto dominance is therefore valuable as an ordering and exploration concept, but it is not inherently designed to identify the single group-preference optimum defined by a stakeholder preference model. The decoupling between optimisation and preference aggregation is most costly when preferences are complex: when stakeholders conflict and no solution is best in all stakeholder directions simultaneously.

Constraint handling in CMOEAs typically follows either the constraint-domination principle \cite{deb2000}, which always prefers feasible solutions over infeasible ones, or penalty functions that add a cost proportional to constraint violation. Both approaches can struggle when the feasible region is very small or when constraint violation correlates strongly with objective quality: infeasible solutions near the feasible boundary carry genetic material that penalty mechanisms destroy \cite{Liang2023A}. More recent approaches such as $\epsilon$-constraint relaxation \cite{TakahamaSakai2006} and archive-based feasibility tracking \cite{Zhang2025} address this, but none integrate constraint handling with direct, PFM-based preference aggregation.

\subsection{The Biased Random-Key Genetic Algorithm}
The Biased Random-Key Genetic Algorithm (BRKGA) \cite{goncalves2013brkga} encodes solutions as vectors of continuous random keys in $[0,1)^n$. A problem-specific decoder maps each chromosome to a feasible (or evaluated) solution. In its most powerful applications, constraint satisfaction is handled entirely within this decoder, decoupling it from the evolutionary operators. This architecture makes BRKGA particularly well-suited to problems with hard, structured feasibility requirements, such as routing, scheduling, and assignment problems where feasibility depends on the combinatorial structure of the solution rather than smooth constraint functions. While conventional penalty-function methods often struggle on such problems because the penalty landscape does not faithfully guide the search toward feasible regions, a well-designed BRKGA decoder sidesteps this by mapping the search space directly to feasible solutions, making feasibility a property of the decoding procedure rather than a gradient to follow.

Multi-objective BRKGA variants have been increasingly explored \cite{londe2025brkga_review}, with existing approaches largely relying on weighted-sum scalarisation, Pareto dominance schemes (such as NSGA-II), lexicographic ordering, or indicator-based selection. These aggregation methods serve exploration and scalar trade-off purposes well, but they were not designed for direct group-preference maximisation: weighted-sum BRKGA, for instance, cannot target an arbitrary point in a non-convex objective space, and it does not perform the affine preference aggregation prescribed by PFM \cite{wolfert2026preferencebasedoptimisationgroupdecisionmaking}. To the best of our knowledge, IMAP-BRKGA is the first multi-objective BRKGA variant to embed PFM-based preference aggregation directly as its fitness signal alongside constraint handling.

\subsection{NSGA-II and the Pareto Selection Paradigm}
Non-dominated Sorting Genetic Algorithm II (NSGA-II) \cite{deb2002fast} remains the most widely used CMOEA. Its defining contribution is fast non-dominated sorting combined with crowding-distance tournament selection, which efficiently maintains a diverse Pareto front approximation across generations. The widespread adoption of NSGA-II makes its Pareto selection paradigm the de facto standard against which new approaches are measured.

The same property that makes NSGA-II powerful for front approximation is the source of its structural indirection for preference-optimal decision problems: Pareto selection has no notion of which part of the front is preferred. Post-hoc extraction of a preferred solution from the final Pareto approximation requires either a secondary optimisation step or a preference model applied to the front after search has terminated. In the latter case, the search has not been guided towards the preferred region. It may have explored regions of the front that are entirely irrelevant to the stakeholders' actual preferences, representing a highly inefficient use of computational resources~\cite{larraga2026}. This paper therefore does not present its continuous instantiation as an NSGA-II variant. IMAP-GA-II keeps only standard real-coded genetic operators also used in many NSGA-II implementations, specifically SBX crossover and polynomial mutation, and removes the two components that make NSGA-II a Pareto-front optimiser: fast non-dominated sorting and crowding-distance selection. The search is consequently guided towards the preference optimum from the first generation rather than towards a diverse Pareto-front approximation.

\subsection{Scalarisation, Preference Direction, and Single-Objective Reduction}
A common alternative to Pareto-front approximation is to reduce the multi-objective problem to a single-objective problem by assigning weights to raw objectives or to normalised objective values. This creates an apparently simple scalar fitness signal, but it is only preference-preserving under restrictive conditions. Objective scales must be commensurate, since a weighted sum of quantities in different units is sensitive to the (often arbitrary) normalisation that makes them comparable, and the mapping from weights to obtained trade-offs is neither linear nor able to reach points on non-convex parts of the Pareto front \cite{das1997drawbacks,marler2010weighted}. Furthermore, the aggregation must respect the measurement structure of preference \cite{Barzilai2010}, and improvement in the scalarised score must correspond to improvement in the stakeholder preference model. When preferences are stakeholder-specific, reversed relative to the default minimisation direction, or non-monotone, this correspondence no longer holds. Minimising a weighted sum of raw objectives then searches for the lowest weighted objective value, not necessarily for the point with highest group preference. Other preference-aware evolutionary approaches incorporate preferences through weighted sums \cite{Kaddani2017Weighted}, reference points \cite{Thiele2009A,Branke2016Using}, ordinal classification \cite{Castellanos-Alvarez2021A}, or by reformulating preferences as constraints \cite{Hou2020Reformulating}; these guide search toward preferred regions of the front, but they retain Pareto dominance or scalarisation as the underlying search signal rather than aggregating stakeholder preference functions directly.

\begin{figure}[t]
    \centering
    \includegraphics[width=0.95\linewidth]{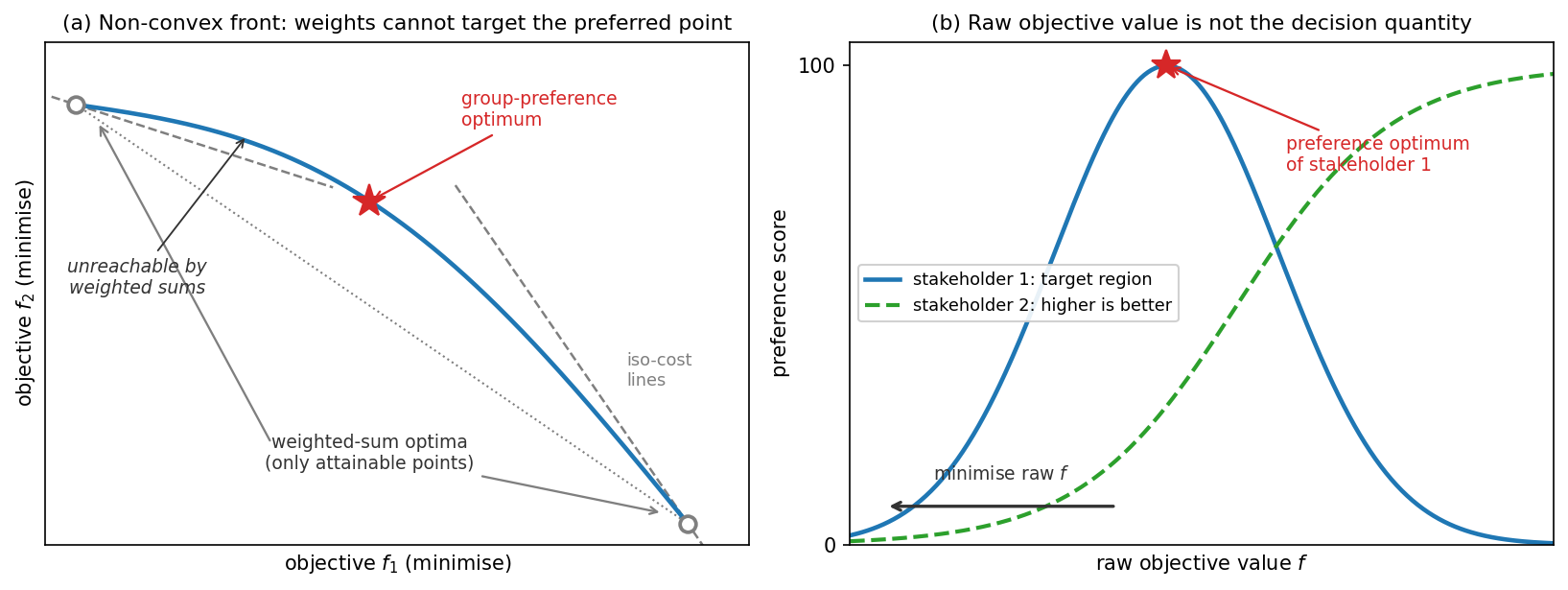}
    \caption{Schematic illustration of why weighted sums of raw objectives are not decision-valid substitutes for preference aggregation. (a)~In objective space, weighted-sum iso-cost lines are supporting lines of the Pareto front: for a non-convex front, every weight vector attains an extreme point, and the concave region containing the group-preference optimum is unreachable for any choice of weights. (b)~In preference space, a stakeholder preference function need not be monotone in the raw objective: a target-region preference peaks at an interior objective value, and a second stakeholder may prefer higher raw values, so minimising the raw objective moves the search away from the preference optimum. The curves are illustrative and not derived from experiment data.}
    \label{fig:scalarisation-geometry}
\end{figure}

The distinction can be seen geometrically (Figure~\ref{fig:scalarisation-geometry}). If the stakeholder preference is identical to minimising a physical performance value, scalarisation may be adequate. But when the front is non-convex, no weight vector can attain its concave region, in which the group-preference optimum may lie (Figure~\ref{fig:scalarisation-geometry}a); and if the stakeholder prefers a target region, or a value that is high for one stakeholder and low for another, then the raw objective value is not itself the decision quantity (Figure~\ref{fig:scalarisation-geometry}b). IMAP-IGS therefore applies stakeholder preference functions before aggregation: the search is driven by preference scores, not by raw objective values whose desirability depends on context. This is why scalarisation is used in the experiments as a paired control rather than as a decision-valid substitute for preference aggregation.

\subsection{Preference Function Modelling and IMAP}
Preference Function Modelling (PFM) \cite{Barzilai2010} provides the measurement-theoretic foundation for multi-criteria preference representation. Preferences are represented on affine scales, so meaningful aggregation must preserve preference differences rather than operate on arbitrary raw objective magnitudes. Within the ODESYS/IMAP line of work, Van Heukelum et al.~\cite{Van_Heukelum2023Socio-technical} introduced IMAP as part of a socio-technical design optimisation methodology that integrates object capability and stakeholder desirability. Wolfert~\cite{wolfert2026preferencebasedoptimisationgroupdecisionmaking} further frames decision-valid preference-based optimisation through the conditions of Preference-Key, Integration, Association, and Uniqueness.

The specific aggregator used in this paper is the weighted centroid of z-normalised preference scores. Its construction follows Wolfert's unique preference aggregation argument: raw stakeholder preference scores are mapped into a common linear preference space through z-normalisation, after which the weighted centroid provides the affine aggregate \cite{wolfert2026uniqueaggregation}. Under this interpretation, weighted sums of raw objectives are not equivalent to IMAP because objective values are not preference scores and their units need not be commensurate. Pareto-based selection is also structurally different: it provides dominance ordering, but it does not itself perform PFM-consistent aggregation to identify a unique group-preference optimum.

IMAP has been used as both a post-hoc scoring tool and, in the Preferendus tool, as an evolutionary fitness signal. Van Heukelum et al.~\cite{Van_Heukelum2023Socio-technical} demonstrate the latter using a GA with an inter-generational solver in infrastructure design applications. The remaining gap is therefore not the first use of IMAP inside any GA, but the absence of a generally applicable, host-agnostic framework that is tested across structurally diverse algorithms and highly constrained problem classes. The primary technical obstacle to such a generalised embedding is intergenerational non-stationarity: z-normalisation depends on the mean and standard deviation of the current evaluation pool, which shifts every generation as the population evolves. This constant shifting renders direct cross-generation score comparisons mathematically invalid. Section~\ref{sec:methodology} derives the current-best re-insertion mechanism that resolves this, allowing IMAP to function as a robust, host-agnostic fitness signal.

\subsection{Benchmark Problem Domains}

\paragraph{DAS-CMOP.}
The Difficulty Adjustable and Scalable Constrained Multi-Objective Problems (DAS-CMOP) benchmark suite was introduced by Fan et al.\ \cite{Fan_2020_DASCMOP} specifically to expose two structural weaknesses of CMOEAs under hard feasibility constraints: their difficulty in \emph{discovering} feasible regions that are small, disconnected, or far from the unconstrained optima, and their slow \emph{convergence} when infeasible barriers block the path to the Pareto front. The suite comprises nine test problems: DAS-CMOP1--6 are bi-objective and DAS-CMOP7--9 are tri-objective. Each problem is parameterised by a difficulty tuple $(\eta_1, \eta_2, \eta_3, \eta_4)$ that independently controls the infeasibility ratio: the fraction of the objective space that is infeasible and the geometry of the feasible region. This design enables a controlled gradient from nearly unconstrained settings to configurations where the feasible region is reduced to narrow bands or isolated patches, with infeasibility ratios exceeding 95\% at the highest difficulty levels. The explicit goal of this parameterisation is to create constraint landscapes that stress the feasibility-maintenance mechanisms of CMOEAs rather than merely testing convergence speed on easy problems.

Existing work benchmarking CMOEAs on DAS-CMOP includes C-TAEA \cite{li2019ctaea}, MOEA/D-CDP \cite{Fan2016moea/d-cdp}, NSGA-II \cite{deb2002fast}, and NSGA-III \cite{deb2014evolutionary}. All of these methods show degraded performance at high difficulty levels: as the feasible region shrinks, Pareto-based selection continues to explore the full approximated front, including infeasible and near-boundary regions irrelevant to any particular preference structure, rather than concentrating search on the preference-aligned subset of the front. The benchmark has no established preference-based baseline because no prior work embeds a preference signal directly into the CMOEA applied to DAS-CMOP. Post-hoc preference selection from a Pareto front approximation is the implicit but unevaluated default. DAS-CMOP is used here as the continuous constrained benchmark domain because its controllable difficulty and publicly available test instances permit rigorous, reproducible comparison of preference-guided and Pareto-based search across the full range from lightly to severely constrained settings.

\paragraph{MO-VRPTW.} The vehicle routing problem with time windows (MO-VRPTW) is a classical NP-hard combinatorial optimisation problem \cite{solomon1987VRPTW} with extensive literature. Kumar \cite{KUMAR2014MO_VRPTW_FAGA} extends this to a multi-objective problem, defining three objectives: total distance travelled, number of vehicles used and route balance. NSGA-II \cite{deb2002fast}, NSGA-III \cite{deb2014evolutionary} and SPEA2 \cite{zitzler2001spea2} variants represent the state of the art in Pareto-based routing. Most existing methods return Pareto fronts with large numbers of non-dominated route plans. No prior MO-VRPTW work embeds preference-based selection directly into the evolutionary search \cite{moghdani2021}. The routing decision-maker is left to select a single route plan from a large non-dominated set without algorithmic support for the preference aggregation step.

\paragraph{Maritime and offshore scheduling.} The HVASP addressed in this paper belongs to a class of heterogeneous fleet scheduling problems arising in offshore energy, marine construction, and logistics \cite{Sharif_multi_objective_vessel_routing_problems_with_safety, Guo2023Integrated, Homsi2020Industrial}. These problems involve discrete assignment of vehicles or vessels to tasks under strict eligibility and time-window constraints. Hybrid constraint programming and metaheuristic approaches have been applied to related scheduling problems \cite{Perez2023A, Iksan2024Optimizing, Lagaros2023Constraint}, but without PFM-compliant preference aggregation.

\subsection{Research Gap}

Three distinct gaps motivate this work.

First, although Van Heukelum et al.~\cite{Van_Heukelum2023Socio-technical} embedded IMAP in the Preferendus solver, no existing work has formulated and tested IMAP as a host-agnostic fitness signal across structurally different GA architectures. Most widely used CMOEAs rely on Pareto-front exploration or scalarisation, search signals that pursue a different goal from PFM-based group preference aggregation~\cite{wolfert2026preferencebasedoptimisationgroupdecisionmaking}. The methodological gap is therefore the absence of a general framework that embeds PFM-based preference aggregation as the search signal itself while preserving the host algorithm's variation operators.

Second, existing preference-based evolutionary work has not yet shown reliable identification of preference-optimal solutions across the combination targeted here: severe feasibility restrictions, structurally different search architectures, and genuine stakeholder conflict. Van Heukelum et al.~\cite{Van_Heukelum2023Socio-technical} show that IMAP can be embedded in a GA for infrastructure design, but not that such embedding scales across continuous constraint-stress benchmarks, combinatorial routing, and real-world heterogeneous scheduling. Post-hoc application of IMAP to Pareto-front outputs can only select among the points that the Pareto search actually returned. If the evolutionary search was not guided towards the group preference optimum, the returned approximation may under-sample or entirely miss the preference-optimal region.

Third, no existing PFM-compliant method handles the all-infeasible-pool problem through the preference formalism itself. Existing CMOEAs address fully infeasible populations via $\epsilon$-constraint relaxation \cite{TakahamaSakai2006}, constraint-dominance principle \cite{deb2000}, or penalty functions applied to objectives: all mechanisms that operate externally to the aggregation step. None express constraint proximity as a native preference function, leaving the aggregator blind to the gradient toward feasibility precisely when the population is hardest to navigate: in early generations on tight instances, when no feasible solution has yet been found.

\section{Problem Formulation}
\label{sec:problem_formulation}

\subsection{Constrained Multi-objective Optimisation Decision Problems}
\label{subsec:CMODP}

Without loss of generality, CMOPs are described as follows \cite{Hao_CMOPs}:

\begin{equation}
\begin{aligned}
& \text{minimise}   & & F(\mathbf{x}) = [f_1(\mathbf{x}), f_2(\mathbf{x}), \dots, f_o(\mathbf{x})] \\
& \text{subject to} & & g_n(\mathbf{x}) \leq 0, \quad \text{for } n = 1, 2, \dots, N \\
&                   & & h_e(\mathbf{x}) = 0, \quad \text{for } e = 1, 2, \dots, E \\
&                   & & \mathbf{x} = [x_1, x_2, \dots, x_D] \in S,
\end{aligned}
\end{equation}

\noindent where $S$ is the decision space and $\mathbf{x}$ is the $D$-dimensional decision vector, $F(\mathbf{x})$ is the objective vector that contains $o$ objectives, $f_i(\mathbf{x})$ is the $i$-th objective function, $g_n(\mathbf{x})$ and $h_e(\mathbf{x})$ are inequality and equality constraints, and $N$ and $E$ are the constraint counts. (The index $m$ is reserved for the number of objectives of a specific problem instance in the remainder of this paper.)

When using CMOPs to reflect real-world design--decision processes, the decision target is not the minimisation of raw objective values as such, but the identification of the feasible solution most preferred by the group of decision-makers. In benchmark terminology the quantities $f_i$ are commonly called objective functions; in the IMAP/ODESYS framing they are treated as system performance dimensions to which stakeholder preference functions are attached.

Following Barzilai~\cite{Barzilai2010}, each decision-maker defines a preference function on an affine scale, for our purposes $[0, 100]$. This preference function $P$ maps performance values to a preference score, where a higher score is more preferred. The interval is a convention rather than a measurement claim: like the zero point of the Celsius temperature scale, the score $0$ is a chosen anchor, the least preferred outcome under consideration, not an absence of preference. Only \emph{differences} between preference scores therefore carry meaning, which is why any admissible aggregation must be restricted to difference-preserving (affine) operations. Building on ODESYS/IMAP \cite{Van_Heukelum2023Socio-technical,wolfert2026preferencebasedoptimisationgroupdecisionmaking} and the unique preference aggregation formulation \cite{wolfert2026uniqueaggregation}, we define the Constrained Multi-objective Optimisation Decision Problem (CMODP) as:

\begin{equation}
\label{eq:CMODP}
\begin{aligned}
& \text{maximise}  & & \mathbf{A}(\mathbf{x}) = \mathbf{A}(P(F(\mathbf{x})),\, \mathbf{w})\\
& \text{where}     & & F(\mathbf{x}) = [f_1(\mathbf{x}), f_2(\mathbf{x}), \ldots, f_o(\mathbf{x})] \\
&                  & & P(F(\mathbf{x})) = \big[p_{k,i}(F(\mathbf{x}))\big], \quad k = 1, \dots, K, \quad i = 1, \dots, o\\
& \text{subject to}& & g_n(\mathbf{x}) \leq 0, \quad n = 1, 2, \dots, N\\
&                  & & h_e(\mathbf{x}) = 0, \quad e = 1, 2, \dots, E\\
&                  & & \mathbf{x} \in S,\\
&                  & & \sum_{k, i} w_{k, i} = 1,
\end{aligned}
\end{equation}

\noindent where $p_{k,i}$ is the preference function of decision-maker $k$ over performance dimension $i$ (with $K$ decision-makers in total), $\mathbf{w} = [w_{k,i}]$ are the preference weights, and $\mathbf{A}(\cdot)$ is the affine aggregator, defined as the weighted centroid of $z$-normalised preference scores $p_{k, i}(F(\mathbf{x}))$. The flat weight vector $\mathbf{w}$ is derived from a nested structure: each decision-maker $k$ carries an importance weight $W_k$ and distributes intra-stakeholder objective weights $w_{i \mid k}$ over the performance dimensions, so that
\begin{equation}
\label{eq:nested-weights}
    w_{k,i} = W_k \cdot w_{i \mid k}, \qquad \sum_{k=1}^{K} W_k = 1, \qquad \sum_{i=1}^{o} w_{i \mid k} = 1 \;\; \forall k,
\end{equation}
\noindent which implies the flat normalisation $\sum_{k,i} w_{k,i} = 1$ in \eqref{eq:CMODP}. The stakeholder definitions in the remainder of this section state $W_k$ and $w_{i \mid k}$ in this nested form. The aggregator is:

\begin{equation}
\label{eq:a-fine-aggregator}
    \mathbf{A}(\mathbf{x}) = \sum_{k, i}w_{k, i} \cdot z_{k, i}, \qquad
    z_{k, i} = \frac{p_{k, i} - \mu_{k, i}}{\sigma_{k, i}},
\end{equation}

\noindent where $p_{k,i}$ is the preference score assigned by decision-maker $k$ to performance dimension $i$, and $\mu_{k,i}$ and $\sigma_{k,i}$ are the mean and standard deviation of that decision-maker's preference scores for dimension $i$ across the evaluation pool \cite{wolfert2026uniqueaggregation}. The normalisation statistics are computed per $(k,i)$ pair, never pooled across decision-makers: each $p_{k,i}$ is measured on its own affine scale, and pooling scores from different stakeholders before normalising would distort the individual scales prior to aggregation. In the implementation, $\sigma_{k,i}$ is regularised as $\max(\sigma_{k,i}, \delta)$ with $\delta = 10^{-6}$. This guard has a clean interpretation rather than being an arbitrary patch: $\sigma_{k,i} = 0$ occurs precisely when every candidate in the pool receives the same score on that dimension, in which case the numerator $p_{k,i} - \mu_{k,i}$ is also zero and the guarded z-score evaluates to exactly $0$. A zero-variance dimension carries no ranking information, and it contributes nothing to the aggregate; a fully degenerate pool, in which all candidates tie on every dimension, yields a constant aggregate for all candidates. The normalisation therefore degrades gracefully as the population converges. The unique optimal solution is:

\begin{equation}
\label{eq:imap_star}
    \mathbf{x}^\ast = \arg\max_{\mathbf{x} \in S_f}\,
    \mathbf{A}(\mathbf{x}),
\end{equation}

\noindent where $S_f \subseteq S$ is the feasible decision space, defined formally as
\begin{equation}
\label{eq:feasible-set}
    S_f = \left\{ \mathbf{x} \in S \;\middle|\; g_n(\mathbf{x}) \leq 0 \;\; \forall n \in \{1, \dots, N\} \;\text{ and }\; h_e(\mathbf{x}) = 0 \;\; \forall e \in \{1, \dots, E\} \right\}.
\end{equation}
\noindent The affine aggregation in~\eqref{eq:a-fine-aggregator} is the decision criterion used throughout the paper: $\mathbf{x}^\ast$ is the feasible design--decision point with maximum aggregated group preference.

This formulation fixes the decision criterion before any algorithm is run. Once the stakeholder preference functions and weights are part of the CMODP, all algorithms can be evaluated by the same question: which feasible decision returned by the algorithm has the largest value of $\mathbf{A}$? Algorithms may differ in their search policy: one may directly maximise $\mathbf{A}$, another may approximate a Pareto front, and another may optimise a weighted scalar objective. These are alternative solution strategies for the same design--decision problem. Consequently, evaluating returned representatives by IMAP is not an algorithm-specific advantage; it is the evaluation criterion induced by the problem formulation itself.

\subsection{DAS-CMOP}
\label{subsec:das-cmop}
The DAS-CMOP suite \cite{Fan_2020_DASCMOP} provides a scalable test framework in which constraint difficulty is tunable via three independent parameters. The $D$-dimensional decision vector is partitioned as $\mathbf{x} = [\mathbf{x}_{1:m-1},\, \mathbf{x}_{m:D}] \in S$, where $\mathbf{x}_{1:m-1}$ governs the shape of the Pareto front and $\mathbf{x}_{m:D}$ governs distance from it. Each objective takes the form

\begin{equation}
    f_i(\mathbf{x}) = \alpha_i(\mathbf{x}_{1:m-1}) + \beta_i(\mathbf{x}_{1:m-1}, \mathbf{x}_{m:D}), \quad i = 1, 2, \ldots, m,
\end{equation}

\noindent where $\alpha_i$ is a shape function and $\beta_i \geq 0$ is a non-negative distance function, so the unconstrained Pareto front is attained at $\beta_i = 0$ for all $i$.

Three constraint types parameterise the difficulty of the suite. Type-I constraints ($K$ in total) restrict the sub-vector $\mathbf{x}_{1:m-1}$, fragmenting the Pareto front into disconnected segments. Type-II constraints ($P$ in total) bound the feasible range of the distance functions $\beta_i$, directly shrinking the feasible fraction of the search space. Type-III constraints ($Q$ in total) impose infeasible obstacles on the objective values themselves, impeding convergence towards the Pareto front.

Following the CMODP formulation in \eqref{eq:CMODP}, we extend DAS-CMOP by introducing preference functions $p_{k,i}$ for each decision-maker $k$:

\begin{equation}
\begin{aligned}
& \text{maximise}   & & \mathbf{A}(\mathbf{x}) = \mathbf{A}\!\left(P(F(\mathbf{x})), \mathbf{w}\right) \\
& \text{where}      & & F(\mathbf{x}) = [f_1(\mathbf{x}),\ f_2(\mathbf{x}), \ldots,\ f_m(\mathbf{x})] \\
& \text{subject to} & & g_n^{(I)}(\mathbf{x}) \leq 0, \quad n = 1,\ 2, \ldots, K \\
&                   & & g_n^{(II)}(\mathbf{x}) \leq 0, \quad n = 1,\ 2, \ldots, P \\
&                   & & g_n^{(III)}(\mathbf{x}) \leq 0, \quad n = 1,\ 2, \ldots, Q \\
&                   & & \mathbf{x} \in S, \quad \sum_{k,i} w_{k,i} = 1,
\end{aligned}
\end{equation}

\noindent where the affine aggregator $\mathbf{A}(\cdot)$ and the z-normalised preference scores $z_{k,i}$ are as defined in \eqref{eq:a-fine-aggregator}.

\paragraph{Preference functions}
For each objective $f_i$, analytical bounds $[\ell_i, u_i]$ are derived from the DAS-CMOP structure. With $D = 30$ variables and $m$ objectives, the distance sub-vector contains $\kappa = D - (m-1)$ components. For DAS-CMOP1--3 and DAS-CMOP9, the distance function satisfies $\beta(\mathbf{x}) \in [0, \kappa]$, giving $\ell_i = 0$ and $u_i = 1 + \kappa$. For DAS-CMOP4--8, the distance function satisfies $\beta(\mathbf{x}) \in [0, 2.25\kappa]$, giving $\ell_i = 0$ and $u_i = 1 + 2.25\kappa$. These bounds are independent of the difficulty parameters, which control only the constraints, and are identical under both stakeholder configurations.

Two stakeholders are defined; their importance weights and objective weights are fixed across both configurations for each problem family (Section~\ref{subsec:stakeholder-configs}), and only the preference-function direction changes between Configuration C and Configuration A.
For bi-objective problems (DAS-CMOP1--6):
\begin{itemize}
    \item $S_1$ (importance weight $W_1 = 0.3$): equal objective weights $w_{i \mid 1} = [0.5, 0.5]$.
    \item $S_2$ (importance weight $W_2 = 0.7$): equal objective weights $w_{i \mid 2} = [0.5, 0.5]$.
\end{itemize}
For tri-objective problems (DAS-CMOP7--9):
\begin{itemize}
    \item $S_1$ (importance weight $W_1 = 0.6$): equal objective weights $w_{i \mid 1} = [1/3, 1/3, 1/3]$.
    \item $S_2$ (importance weight $W_2 = 0.4$): unequal objective weights $w_{i \mid 2} = [0.2, 0.3, 0.5]$ emphasising $f_3$.
\end{itemize}

Preference-function directions (minimising vs.\ maximising) depend on the stakeholder configuration and are defined in Section~\ref{subsec:stakeholder-configs}.

\subsection{Multi-Objective Vehicle Routing Problem with Time Windows}
\label{subsec:vrptw}
The MO-VRPTW~\cite{KUMAR2014MO_VRPTW_FAGA} models a central depot dispatching $K$ homogeneous vehicles to serve $N$ customers, each with a known demand $q_i$, service time $s_i$, and time window $[e_i, l_i]$. The decision vector $\mathbf{x}$ encodes the assignment of customers to routes and the sequence within each route. Travel distances and times between nodes are computed as Euclidean distances.

Three objectives are minimised:
\begin{subequations}
    \begin{align}
        & f_1(\mathbf{x}) = \sum_{j=1}^{K} d_j  & & \text{(total distance)}\\
        & f_2(\mathbf{x}) = K                   & & \text{(number of vehicles)}\\
        & f_3(\mathbf{x}) = d_{\max} - d_{\min} & &  \text{(route balance)},
    \end{align}
\end{subequations}

\noindent where $d_j$ is the total distance of route $j$, and $d_{\max}$, $d_{\min}$ are the longest and shortest route distances respectively. Mapping to the CMODP in~\eqref{eq:CMODP}:

\begin{equation}
\begin{aligned}
& \text{maximise}   & & \mathbf{A}(\mathbf{x}) = \mathbf{A}\!\left(P(F(\mathbf{x})), \mathbf{w}\right), & & \\
& \text{where}      & &  F(\mathbf{x}) = [f_1(\mathbf{x}),\ f_2(\mathbf{x}),\ f_3(\mathbf{x})],         & & \\
& \text{subject to} & &  \sum_{i \in \text{route } j} q_i \leq C, \quad \forall j = 1,\ 2, \ldots, K        & & \text{(capacity)} \\
&                   & & e_i \leq a_i + w_i \leq l_i, \quad \forall i                                    & & \text{(time windows)} \\
&                   & & \sum_{i \in \text{route } j} (t_i + w_i + s_i) \leq t_{\max}, \quad \forall j   & & \text{(maximum route time)} \\ 
&                   & & a_0 = w_0 = s_0 = q_0 = 0,                                                      & & \text{(depot)}\\
& & & \sum_{k,i} w_{k,i} = 1 & &
\end{aligned}
\end{equation}

\noindent Here $C$ is the vehicle capacity, $a_i$, $w_i$, and $s_i$ denote the arrival time, waiting time, and service time at customer $i$, $t_i$ the travel time to customer $i$, and $t_{\max}$ the maximum route duration; index $0$ denotes the depot.

\paragraph{Preference functions}
Objective bounds $[\ell_i, u_i]$ are computed analytically per Solomon instance without evaluating any solutions. For total distance $f_1$: $\ell_1$ is the minimum spanning tree (MST) weight over the depot and all customers (Prim's algorithm), and $u_1 = 2\sum_{c} d(\text{depot}, c)$, which is always reachable by sending one vehicle per customer. For vehicle count $f_2$: $\ell_2 = \lceil \sum_c q_c / C \rceil$ (bin-packing lower bound) and $u_2 = N$. For route balance $f_3 = d_{\max} - d_{\min}$: $\ell_3 = 0$ (achievable when all routes have equal length) and $u_3 = 2 \sum_{j \in \mathcal{J}} d(\text{depot}, c_j)$, where $\mathcal{J}$ contains the $n_{\max}$ farthest customers fitting on a single route.

Two stakeholders are defined:
\begin{itemize}
    \item $S_1$ (importance weight $W_1 = 0.6$): equal objective weights $w_{i \mid 1} = [1/3, 1/3, 1/3]$.
    \item $S_2$ (importance weight $W_2 = 0.4$): unequal objective weights $w_{i \mid 2} = [0.2, 0.3, 0.5]$ placing greater emphasis on route balance $f_3$.
\end{itemize}

Preference-function directions (minimising vs.\ maximising) depend on the stakeholder configuration and are defined in Section~\ref{subsec:stakeholder-configs}.

\subsection{HVASP}
The Heterogeneous Vessel Allocation and Scheduling Problem (HVASP) is a real-world CMOP in which a fleet of $|\mathcal{V}|$ vessels must be assigned and sequenced across a portfolio of $|\mathcal{A}|$ activities subject to timing, precedence, and vessel-capability constraints. The decision vector $\mathbf{x} = (x_1, \ldots, x_6)$ encodes activity start times, maintenance location choices, vessel-to-role assignments, sequencing decisions, sequence-start indicators, and inter-role sailing speeds. The exogenous parameter vector $\mathbf{y} = (y_1, \ldots, y_{15})$ captures activity durations, time windows, spatial data, precedence relations, vessel capabilities, and fuel profiles; full definitions are provided in Appendix~\ref{appendix:HVASP}.

Two objectives are minimised:

\begin{subequations}
    \begin{align}
        & f_1(\mathbf{x}) = \sum_{r \in \mathcal{R}} \sum_{r' \in \mathcal{R}} x_4[r,r'] \cdot \gamma\!\left(x_3[r], r, r', x_6[r,r'], \delta_{r,r'}\right) & &  \text{(total cost)},\\
        & f_2(\mathbf{x}) = \max_{a \in \mathcal{A}}\!\left(x_1[a] + y_1[a]\right) - \min_{a \in \mathcal{A}}\, x_1[a] & &  \text{(makespan)},
    \end{align}
\end{subequations}

\noindent where the transition cost function $\gamma$ is defined in Appendix~\ref{appendix:HVASP}. The problem is subject to 13 constraints $g^{(m')}_f(\mathbf{x},\mathbf{y}) \leq 0$, $m'=1,\ldots,13$, covering activity precedence and vessel exclusivity ($g^{(1)}$-$g^{(2)}$), role sequencing structure ($g^{(3)}$–$g^{(9)}$), and path integrity per vessel ($g^{(10)}$–$g^{(13)}$), as detailed in Appendix~\ref{appendix:HVASP}. The feasible space is $S_f = \{(\mathbf{x},\mathbf{y}) \mid g^{(m')}_f(\mathbf{x},\mathbf{y}) \leq 0,\ m'=1,\ldots,13\}$.

Mapping to the CMODP in~\eqref{eq:CMODP}, the HVASP is stated as

\begin{equation}
\begin{aligned}
    & \text{maximise}   & & \mathbf{A}(\mathbf{x}) = \mathbf{A}\!\left(P(F(\mathbf{x})),\,\mathbf{w}\right) \\
    & \text{where}      & & F(\mathbf{x}) = [f_1(\mathbf{x}), f_2(\mathbf{x})] \\
    & \text{subject to} & & \mathbf{x} \in S_f\\
    &                   & & \sum_{k,i} w_{k,i} = 1, \\
\end{aligned}
\end{equation}

\noindent with affine aggregator $\mathbf{A}(\cdot)$ and z-normalised scores $z_{k,i}$ as defined in~\eqref{eq:a-fine-aggregator}. The parameter vector $\mathbf{y}$ is fixed per instance, so feasibility of $\mathbf{x}$ is understood with respect to the given $\mathbf{y}$.

\paragraph{Preference functions}
Objective bounds $[\ell_i, u_i]$ for cost and makespan are computed exactly by running the single-objective constraint solver (CP Phase~A: minimise and maximise each objective independently). Because the bounds are instance-specific, derived from the synthetic instance's parameters rather than from closed-form expressions, this CP pre-solve is required for each instance before any metaheuristic run. Stakeholder configurations and preference-function directions are defined in Section~\ref{subsec:stakeholder-configs}.

\section{Methodology}
\label{sec:methodology}

This section develops the IMAP-IGS framework in three movements. Sections~\ref{subsec:preference_function_construction} and~\ref{subsec:preference_validity} construct the stakeholder preference functions and delimit what the preference model does and does not claim. Sections~\ref{subsec:imap-ga_framework} to~\ref{subsec:cv-pref} present the framework itself: IMAP as a scalar fitness replacement, the intergenerational evaluation problem created by pool-relative scoring and its solution, and the constraint-violation preference function that steers the search toward feasibility. Sections~\ref{subsec:imap-brkga} and~\ref{subsec:imap-gaii} then instantiate the framework in two structurally opposite host algorithms.

\subsection{Preference Function Construction}
\label{subsec:preference_function_construction}
For each problem instance, IMAP preference functions are constructed per stakeholder
and per objective. Each function $p_{k,i}: [l_i,\; u_i] \to [0,\; 100]$ maps an
objective value to a preference score, where $[l_i,\; u_i]$ are instance-specific
objective bounds. For a stakeholder whose preference direction agrees with the objective's raw minimisation formulation, the function is linearly decreasing:
$p_{k,i}(v) = 100 \cdot \frac{(u_i - v)}{(u_i - l_i)}$. For a stakeholder whose preference direction opposes it, so that the dimension is functionally a maximisation for that stakeholder, the function is reversed: $p_{k,i}(v) = 100 \cdot \frac{(v - l_i)}{(u_i - l_i)}$.

For DAS-CMOP and MO-VRPTW, objective bounds are computed analytically from the
problem structure. For HVASP, bounds are computed instance-specifically by the CP
solver's Phase A (single-objective minimisation and maximisation runs per objective),
ensuring preference functions are calibrated to the actual achievable range for each
instance.

\subsection{Preference-Function Validity and Scope}
\label{subsec:preference_validity}
The IMAP-IGS framework assumes that stakeholder preference functions and stakeholder weights are accepted as part of the decision model. IMAP-IGS therefore answers a conditional question: given these preference functions, which feasible solution has the highest aggregated group preference? It does not, by itself, validate whether the elicited functions perfectly represent stakeholder psychology, organisational priorities, or future changes in those priorities.

The linear functions used in the experiments are controlled preference models chosen to isolate the effect of embedding preference aggregation into search. They should not be read as a claim that all real decision-makers have linear preferences. In a deployed decision-support setting, the same framework could use elicited piecewise-linear, concave, convex, target-seeking, or otherwise stakeholder-specific preference functions, provided that they are expressed on the affine preference scale required by PFM. Robustness to misspecified functions can then be studied by perturbing slopes, breakpoints, target values, and stakeholder weights.

In real decision processes, preference elicitation is itself part of the optimisation cycle rather than a one-off pre-processing step. If the highest-preference solution returned by IMAP-IGS is not satisfactory to the decision-maker, this does not necessarily indicate an algorithmic failure. It may reveal that the stated preference functions or weights did not accurately capture the decision-maker's true priorities. The appropriate response is then to revisit the elicitation step, reshape the relevant preference functions, and re-run the search with the updated decision model. This iterative loop is consistent with the socio-technical design process described in Odesys and Preferendus \cite{Van_Heukelum2023Socio-technical, Wolfert2023}: the optimiser provides feedback about the consequences of the stated preferences, and that feedback can help stakeholders refine what they actually mean by a preferred solution.

Appendix~\ref{appendix:pareto-imap} clarifies the relationship between IMAP optimality and Pareto efficiency. In the special case of aligned monotone preference functions, a complete Pareto front followed by IMAP selection can in principle recover the same decision point, provided the final selection is scored in a shared evaluation pool. Under conflicting, reversed, or non-monotone preference functions, this equivalence no longer holds in general, which is precisely the regime targeted by IMAP-IGS.

\subsection{The IMAP-IGS Framework as a Scalar Fitness Replacement}
\label{subsec:imap-ga_framework}

The IMAP-IGS framework adapts any population-based evolutionary algorithm that relies on scalar fitness evaluation. Candidate solutions are generated by the host algorithm, decoded into performance values, mapped to stakeholder preference scores, augmented with a constraint-violation preference score for search-time feasibility steering, z-normalised within the current evaluation pool, and aggregated through the IMAP weighted-centroid operator. The resulting IMAP score replaces the host algorithm's fitness signal for selection. All evolutionary operators, such as chromosome representation, crossover, mutation, repair, and decoding, are inherited unchanged from the host algorithm.

Throughout this section, $A(\cdot)$ denotes the affine aggregation operator defined in Eq.~\eqref{eq:a-fine-aggregator}: the weighted centroid of the z-normalised stakeholder preference scores. A critical property of this aggregator is that its scores are \emph{pool-relative}: $A(\mathbf{x})$ is computed via z-normalisation over the current evaluation pool. Thus, score magnitudes are not comparable across different pools. IMAP scores are meaningful only as a ranking device within the pool in which they were computed. This has a direct consequence for the framework design: selection must be based on within-pool ranking, and the best currently known feasible solution $\mathbf{x}^\dagger$ must be re-inserted into the pool before scoring. The \emph{evaluation pool} of generation $t$ is thus $\mathcal{P}^{\mathrm{ev}}_t = \mathcal{P}_{t-1} \cup \{\mathbf{x}^\dagger\}$: the current population together with the re-inserted current best, and all IMAP scores of generation $t$ are computed over $\mathcal{P}^{\mathrm{ev}}_t$. This current-best re-insertion ensures that $\mathbf{x}^\dagger$ is evaluated in the same normalisation context as the live population. Infeasible individuals are ranked below all feasible ones whenever feasible candidates are present. When the pool is entirely infeasible, the constraint-violation preference function $p_{cv}$ (Section~\ref{subsec:cv-pref}) ensures the ranking remains discriminative even in the absence of a feasible reference.

Algorithm~\ref{alg:imap-ga} presents the general IMAP-IGS pseudocode. The host-GA operators are encapsulated in a single \textsc{EvolvePopulation}($\mathcal{P}_{t-1}$, $A(\cdot)$) call. $A(\cdot)$ replaces whatever scalar fitness signal the host would otherwise use for selection.

The implementation used in the experiments applies a two-pass ranking step to improve convergence. First, all candidates in the evaluation pool are scored and ranked by IMAP and their generation-local scores are rescaled to $[0,100]$. Candidates whose rescaled score is at most a threshold $\tau=40$ are marked for exclusion from selection and from the current-best update. This discard is deliberately guarded so that it cannot deplete the population. The threshold is applied only when enough candidates clear it: if the discard would leave fewer candidates than the nominal population size, it is skipped entirely for that generation and the full pool proceeds to selection. When feasible candidates are scarce, they are additionally exempted from the threshold, and the re-inserted current-best solution, when present, is always retained so that the algorithm never discards its best currently known feasible representative solely because the generation-local normalisation context changed. Note that $\tau$ is a threshold on the rescaled aggregate score, not a population quantile: the fraction of candidates below it depends on the score distribution of the pool and varies per generation, and under the population-size guard the filter engages only when the pool is of sufficient quality to afford it. Second, IMAP scores are recomputed over the retained candidates, together with the re-inserted current-best solution when present, and this second ranking is used for selection. Because IMAP scores are pool-relative, the re-ranking step is important: after low-preference candidates are removed, the normalisation context changes and the retained candidates must be compared within the filtered pool. Population diversity is not governed by this filter alone: the host algorithm's exploration operators act on the full population every generation regardless (in IMAP-BRKGA a fixed fraction of fresh random mutants is injected each generation, and in IMAP-GA-II standard SBX crossover and polynomial mutation continue to generate variation), so the filter sharpens selection pressure without truncating the source of new genetic material. The threshold is an algorithmic convergence device, not part of the final evaluation metric; all reported $Z^*$ values are computed post-hoc under the common evaluation protocol in Section~\ref{subsec:metrics}. The sensitivity of the framework to $\tau$ is discussed together with the other algorithmic constants at the end of Section~\ref{subsec:cv-pref}.

\begin{algorithm}[t]
\caption{IMAP-IGS (framework)}
\label{alg:imap-ga}
\begin{algorithmic}[1]
\Require Preference functions $\{p_{k,i}\}$, substantive weights $\{w_{k,i}\}$, feasibility weight $w_\mathrm{cv}$, initial population $\mathcal{P}_0$, max generations $T$, threshold $\tau$, host-GA operator \textsc{EvolvePopulation}
\Ensure Best feasible solution $\mathbf{x}^\dagger$
\State $\mathbf{x}^\dagger \leftarrow \varnothing$
       \Comment{current-best re-insertion representative, stored by objective values $F(\mathbf{x}^\dagger)$}
\For{$t = 1, \ldots, T$}
  \State Evaluate $F(\mathbf{x})$ for all $\mathbf{x} \in \mathcal{P}_{t-1}$
  \If{$\mathbf{x}^\dagger \neq \varnothing$}
    \State Re-insert $\mathbf{x}^\dagger$ into the evaluation pool $\mathcal{P}^{\mathrm{ev}}_t = \mathcal{P}_{t-1} \cup \{\mathbf{x}^\dagger\}$
  \EndIf
  \State Append $p_\mathrm{cv}(\mathbf{x})$ to the search-time preference vector for each $\mathbf{x} \in \mathcal{P}^{\mathrm{ev}}_t$
  \State Compute initial IMAP scores over $\mathcal{P}^{\mathrm{ev}}_t$; rank infeasible individuals last when feasible candidates are present
  \State Rescale generation-local scores to $[0,100]$ and mark candidates with score $\leq\tau$ for discard
  \State Skip the discard if it would leave fewer candidates than the population size; exempt scarce feasible candidates; always retain $\mathbf{x}^\dagger$ when present
  \State Recompute IMAP scores over the retained pool
  \State Update $\mathbf{x}^\dagger$ only if a feasible retained individual outscores it within this pool
  \State $\mathcal{P}_t \leftarrow \textsc{EvolvePopulation}(\mathcal{P}_{t-1},\ A(\cdot))$
         \Comment{host-GA operators; second-pass $A(\cdot)$ replaces fitness}
\EndFor
\State Re-insert $\mathbf{x}^\dagger$ into final pool; recompute IMAP scores over
       $\mathcal{P}_T \cup \{\mathbf{x}^\dagger\}$
\State \Return $\mathbf{x}^\dagger$
\end{algorithmic}
\end{algorithm}

\subsection{Intergenerational IMAP Evaluation}
\label{subsec:intergenerational}
The central technical challenge in applying IMAP inside an evolutionary loop is intergenerational non-stationarity. Because z-normalisation depends on $\mu_{k,i}$ and $\sigma_{k,i}$ computed over the current evaluation pool, these parameters shift every generation as the population evolves. A solution scoring 85 in generation $t$ may score 40 in generation $t+1$ solely because the pool composition changed, even if its objective values are identical. Cross-generation comparisons of IMAP scores are therefore invalid, and naive application of IMAP as a fitness function, without accounting for this, does not reliably converge to the preference optimum. This is a property of z-normalisation itself, not of any particular GA architecture; it arises identically whether the host GA is BRKGA, NSGA-II, or any other evolutionary method.

The solution is current-best re-insertion: maintain a current-best solution $\mathbf{x}^\dagger$ across generations, tracked by objective values $F(\mathbf{x}^\dagger)$ rather than IMAP score. Storing $\mathbf{x}^\dagger$ by objective values rather than by score sidesteps the non-stationarity: objective values do not change between generations, so $\mathbf{x}^\dagger$ remains a stable representative regardless of how the normalisation context shifts. At each generation, $\mathbf{x}^\dagger$ is re-inserted into the evaluation pool \emph{before} IMAP scores are computed, ensuring that its score is always evaluated in the same normalisation context as the live population. The current-best re-insertion representative is updated to any new feasible individual that outscores $\mathbf{x}^\dagger$ \emph{within that generation's pool}. The within-pool qualification is essential, since it is the only context in which IMAP score comparisons are valid.

This solution is derived once here and applied identically in both IMAP-BRKGA and IMAP-GA-II. It requires no modification to any evolutionary operator and imposes no assumption about chromosome representation or decoding strategy.

\subsection{Constraint-Violation Preference Function}
\label{subsec:cv-pref}

The framework handles constraints in the preference layer rather than by relying on objective-space penalties as the primary steering mechanism. Constraint-domination ensures that feasible individuals are ranked above infeasible individuals whenever both are present in the evaluation pool. A structural failure arises, however, when \emph{every} individual in the pool is infeasible. In this case constraint-domination has no discriminative effect because there are no feasible candidates to dominate the infeasible ones. If IMAP were computed only from objective-based preference functions, the search would receive no direct information about proximity to the feasible region. On tight instances, particularly during early generations, this failure mode is not hypothetical; it is a regular condition that the framework must handle.

The solution is to define constraint proximity as an auxiliary preference function $p_\mathrm{cv}$ and include it in the search-time preference aggregate with a large feasibility weight. In the experiments, the constraint-violation preference is assigned a pseudo-stakeholder weight of:

\begin{equation}
    w_\mathrm{cv} = 0.50.
    \label{eq:w-cv}
\end{equation}

\noindent The substantive stakeholder preference weights are then rescaled by the remaining mass:

\begin{equation}
    \tilde{w}_{k,i} = (1-w_\mathrm{cv})\,w_{k,i},
    \qquad
    \sum_{k,i} \tilde{w}_{k,i} + w_\mathrm{cv}=1.
    \label{eq:cv-weight-rescaling}
\end{equation}

\noindent The search-time aggregate is therefore:

\begin{equation}
    A_\mathrm{search}(\mathbf{x}) =
    \sum_{k,i} \tilde{w}_{k,i} z_{k,i}(\mathbf{x})
    + w_\mathrm{cv} z_\mathrm{cv}(\mathbf{x}).
    \label{eq:search-aggregate-cv}
\end{equation}

The constraint-violation preference enters the aggregate through the same z-normalisation as the substantive preference dimensions: $z_\mathrm{cv}$ in \eqref{eq:search-aggregate-cv} is the z-score of $p_\mathrm{cv}$ over the current evaluation pool. This is a deliberate design choice with a useful consequence: it makes the effective feasibility pressure \emph{self-adaptive} with respect to the pool's feasibility composition. When feasible solutions are rare, the $p_\mathrm{cv}$ column has high variance and the few feasible candidates receive large positive $z_\mathrm{cv}$ scores, concentrating selection pressure on feasibility exactly when feasibility is the binding bottleneck. As the pool fills with feasible candidates, the column's variance shrinks and the pressure relaxes in favour of the substantive preferences; when the pool is entirely feasible, $\sigma_\mathrm{cv} = 0$ and the regularised z-score (Section~\ref{sec:problem_formulation}) evaluates to zero, so the feasibility dimension deactivates itself and ranking is decided purely by stakeholder preference. The magnitude of this pressure is bounded, not explosive: in the bimodal worst case of a pool of size $n_\mathrm{ev}$ in which $q$ feasible candidates score $100$ and all others share the Zone-2 plateau score, the z-score of the feasible group is exactly $\sqrt{(n_\mathrm{ev}-q)/q}$, and in general no z-score in a pool of size $n_\mathrm{ev}$ can exceed $(n_\mathrm{ev}-1)/\sqrt{n_\mathrm{ev}}$. The behaviour that adaptive penalty schemes construct explicitly, penalty pressure that responds to the current feasible ratio of the population \cite{Lagaros2023Constraint,Rahimi2022A}, is here obtained directly from the normalisation, without a hand-designed schedule. An alternative design would add $p_\mathrm{cv}$ as a raw additive penalty outside the normalisation; this is rejected on the same grounds on which the framework rejects raw-objective scalarisation, as it would mix a technical quantity of arbitrary scale into the affine preference space and reintroduce the commensurability problem that the preference layer exists to avoid.

The value $w_\mathrm{cv} = 0.50$ was chosen on the basis of preliminary experimentation rather than theoretical analysis: at this weight the $p_\mathrm{cv}$ term was found to dominate the search-time aggregate whenever a feasible and an infeasible candidate are compared, so that a feasible solution is almost always preferred over an infeasible one, because feasible candidates receive $p_\mathrm{cv}=100$ while infeasible candidates receive a sharply lower score. A static $w_\mathrm{cv}$ does not, however, imply static feasibility pressure: the effective pressure is the product $w_\mathrm{cv} \cdot z_\mathrm{cv}(\mathbf{x})$, and as described above the $z_\mathrm{cv}$ factor scales automatically with the scarcity of feasible candidates in the pool, vanishing entirely once the pool is feasible throughout. The weight fixes how much of the aggregate the feasibility signal may claim when it is informative, and the normalisation decides, per generation, how informative it is. At the same time, infeasible candidates are not discarded automatically; when no feasible solution is available, $p_\mathrm{cv}$ still supplies a gradient toward feasibility. The constraint-violation preference is not a substantive stakeholder objective: it is an algorithmic feasibility-steering preference used during search. Final reported $Z^*$ scores are computed over feasible representatives using only the stakeholder preference functions defined for the decision problem, so $p_\mathrm{cv}$ is not included in the final reporting aggregate.

The specific instantiation used in this work has a three-zone step-decay shape. Let $\varepsilon > 0$ be a small threshold, $k > 0$ a scale factor, and $p > 0$ a shape exponent:
\begin{equation}
  p_\mathrm{cv}(\mathbf{x}) =
  \begin{cases}
    100 & \text{if } \bar{v} = 0, \\
    10  & \text{if } 0 < \bar{v} \leq \varepsilon, \\[4pt]
    \dfrac{10}{1 + k\,(\bar{v} - \varepsilon)^{p}} & \text{if } \bar{v} > \varepsilon,
  \end{cases}
  \label{eq:cv-pref}
\end{equation}
where $\bar{v}$ is the sum of per-constraint violations, normalised to $[0,1]$.
Zone 1 ($\bar{v} = 0$) maps feasible solutions to the maximum score.
Zone 2 ($0 < \bar{v} \leq \varepsilon$) enforces the required discontinuity: the instant a solution becomes infeasible its score drops to 10, creating an unambiguous gap relative to any feasible solution.
Zone 3 ($\bar{v} > \varepsilon$) provides a strictly positive, monotone-decreasing tail anchored at 10 that decays toward, but never reaches, zero, maintaining a gradient toward feasibility throughout the search.

The discontinuity between Zones 1 and 2 is a design requirement, not an artefact. Feasibility is binary in the underlying decision problem: an infeasible solution is not a valid decision, however small its violation, so no amount of substantive preference may trade off against it. A continuous ramp connecting near-feasible scores to the feasible score would permit exactly that trade, allowing a slightly infeasible candidate to outrank a feasible one in the aggregate, which is the situation constraint-domination exists to prevent. The jump therefore severs the gradient only \emph{across} the feasibility boundary, where no gradient is admissible by construction; \emph{within} the infeasible region, the only region in which the search requires steering toward feasibility, $p_\mathrm{cv}$ is strictly monotone in $\bar{v}$ through the Zone-3 tail. A solution with violation $10^{-7}$ scoring 10 while a feasible solution scores 100 is thus the intended behaviour: the former is not yet a decision candidate, and its path to the boundary is still fully graded. The threshold $\varepsilon = 10^{-6}$ should accordingly be read as a numerical feasibility tolerance rather than a tuned shape parameter: violations at or below floating-point and solver precision are treated as boundary cases and share the Zone-2 plateau, preventing numerical noise from creating spurious ranking distinctions among effectively-boundary candidates. The tail parameters $k = 2.0$ and $p = 1.5$ were derived via preliminary sensitivity analysis to balance two competing algorithmic needs: ensuring that the penalty gradient remains computationally distinct from objective-space noise, while simultaneously avoiding a gradient so steep that it forces premature convergence to sub-optimal local minima on the feasibility boundary.

Because $p_\mathrm{cv}$ acts directly in the preference dimension and receives a large search-time weight, it is the primary feasibility-steering mechanism used by IMAP-IGS. Problem-specific decoders may still repair infeasible chromosomes, reject impossible assignments, or report total constraint violation, but they do not need to translate infeasibility into an artificial objective penalty. This keeps feasibility information inside the same preference formalism as the stakeholder objectives and avoids mixing a technical penalty with the substantive objective values. If an evaluation pool contains only feasible candidates, $p_\mathrm{cv}$ has no discriminative effect; in that case its term is omitted from the within-pool ranking or equivalently treated as a constant contribution.

\paragraph{Fixed algorithmic constants.} Four constants are held fixed across every experiment in this paper: the truncation threshold $\tau = 40$ (Section~\ref{subsec:imap-ga_framework}), the feasibility weight $w_\mathrm{cv} = 0.50$, and the tail parameters $k = 2.0$ and $p = 1.5$ of \eqref{eq:cv-pref}. ($\varepsilon$ is a numerical tolerance, not a tuned constant.) All four were set through preliminary experimentation, are used unchanged across all benchmark domains and both host instantiations, and are not claimed to be optimal. A one-at-a-time ablation of all four constants on a DAS-CMOP subset (Appendix~\ref{appendix:ablation}) supports the choices: no neighbouring configuration outperforms the base beyond sampling noise, the framework is insensitive to $\tau$ and $k$ over the tested ranges, and the only significant deviations are degradations at the $w_\mathrm{cv}$ extremes and at a steeper penalty tail, which fail in precisely the understeering and oversteering directions described above. Interaction effects between the constants and their behaviour on the combinatorial domains remain future work.

\subsection{Instantiation I: IMAP-BRKGA}
\label{subsec:imap-brkga}
The IMAP-IGS framework is now instantiated in two host algorithms, chosen deliberately from opposite ends of the GA design space: a decoder-based random-key GA for combinatorial problems with hard feasibility structures (this section), and a real-coded, operator-based GA for continuous constrained problems (Section~\ref{subsec:imap-gaii}). Demonstrating that the identical IMAP machinery drives both, without modification to either host's variation operators, is what substantiates the host-agnosticism claim; a single instantiation could not distinguish a property of the framework from a property of one architecture.

IMAP-BRKGA specialises IMAP-IGS with \textsc{EvolvePopulation} implemented as biased
random-key crossover with elitist inheritance \cite{goncalves2013brkga}.

\paragraph{Chromosome representation.} Each chromosome is a vector
$\mathbf{r} = [r_1, \ldots, r_n]^\top$ with $r_i \in [0,1)$. A problem-specific
decoder $\mathcal{D}: [0,1)^n \to \mathcal{X}$ maps each chromosome to a decoded
solution. Constraint satisfaction is handled entirely within the decoder: feasibility
is a property of the decoded solution, not the chromosome, so the evolutionary
operators never directly encounter constraints.

\paragraph{Population structure and evolution.} The population of size $|\mathcal{P}|$
is partitioned into an elite set ($\rho_e$ individuals with highest
IMAP scores), a mutant set ($\rho_m$ freshly randomised chromosomes),
and a non-elite set (the remainder). In biased crossover, each gene of an offspring
is drawn from the elite parent with probability $\rho_b$ and from the non-elite parent
otherwise, ensuring elites' genetic material propagates efficiently.

\paragraph{Decoder feasibility information.} Constraint handling in IMAP-BRKGA is
split between the problem-specific decoder and the preference-layer feasibility signal.
The decoder constructs a candidate solution from the random keys and records whether
all hard feasibility requirements are satisfied. If a candidate cannot be decoded into a
fully feasible solution, the decoder reports the corresponding normalised violation
$\bar{v}$. This value is then passed to the constraint-violation preference function
$p_\mathrm{cv}$ from Section~\ref{subsec:cv-pref}. Thus infeasible individuals can still
remain in the population and contribute genetic material, but their selection pressure
is governed by the feasibility preference rather than by adding an artificial
penalty to one of the substantive objectives. The constraint-domination principle
additionally ensures that no infeasible individual can become the current-best re-insertion representative
$\mathbf{x}^\dagger$ when feasible candidates exist.

The complete IMAP-BRKGA procedure is given in Algorithm~\ref{alg:imap-brkga}, which
specialises the generic IMAP-IGS framework (Algorithm~\ref{alg:imap-ga}) by
implementing \textsc{EvolvePopulation} as biased random-key crossover with elitist
inheritance. In the elite-selection step, infeasible individuals are ordered behind
all feasible ones in the second-pass ranking, so an infeasible individual can enter
the elite set only when there are fewer feasible candidates than elite slots.

\begin{algorithm}[t]
\caption{IMAP-BRKGA (Instantiation I - specialises Algorithm~\ref{alg:imap-ga})}
\label{alg:imap-brkga}
\begin{algorithmic}[1]
\Require Preference functions $\{p_{k,i}\}$, substantive weights $\{w_{k,i}\}$, feasibility weight $w_\mathrm{cv}$, population size $n_p$, elite fraction $p_e$, mutant fraction $p_m$, crossover bias $\rho_b > 0.5$, threshold $\tau$, max generations $T$
\Ensure Best feasible solution $\mathbf{x}^\dagger$
\State Initialise $\mathcal{P}_0$: draw $n_p$ chromosomes uniformly from $[0,1)^D$
\State $\mathbf{x}^\dagger \leftarrow \varnothing$
       \Comment{current-best re-insertion representative, stored by objective values $F(\mathbf{x}^\dagger)$}
\For{$t = 1, \ldots, T$}
  \State Decode each $c \in \mathcal{P}_{t-1}$; evaluate $F(\mathbf{x}(c))$
  \If{$\mathbf{x}^\dagger \neq \varnothing$}
    \State Re-insert $\mathbf{x}^\dagger$ into evaluation pool
  \EndIf
  \State Append $p_\mathrm{cv}(\mathbf{x})$ to the search-time preference vector for each $\mathbf{x}$ in the pool
         \Comment{Eq.~\ref{eq:cv-pref}}
  \State Compute initial IMAP scores over the full pool; rank infeasible individuals last when feasible candidates are present
  \State Filter out candidates with rescaled preference score $\leq \tau$ (skipped if fewer than $n_p$ would remain) and recompute IMAP scores
  \State $\mathbf{x}^\dagger \leftarrow \arg\max A(\mathbf{x})$ over feasible
         individuals in the filtered pool
  \State $\mathcal{E}_t \leftarrow$ top $\lfloor p_e \cdot n_p \rfloor$ individuals
         by second-pass IMAP score \Comment{infeasible ranked last}
  \State Sample $\lfloor p_m \cdot n_p \rfloor$ mutant chromosomes uniformly
         from $[0,1)^D$
  \State Fill remaining slots via biased crossover: each gene from a random elite
         parent with probability $\rho_b$, otherwise from a random non-elite parent
  \State $\mathcal{P}_t \leftarrow \mathcal{E}_t \cup
         \text{mutants} \cup \text{crossover offspring}$
\EndFor
\State \Return $\mathbf{x}^\dagger$
\end{algorithmic}
\end{algorithm}

\paragraph{MO-VRPTW decoder.} Each gene $r_i \in [0,1)$ is associated with customer $i$ and encodes its construction priority. Six construction strategies from the VRPTW literature are implemented as separate interchangeable decoders~\cite{solomon1987VRPTW}: the \textbf{I1}, \textbf{I2}, and \textbf{I3} cheapest-insertion heuristics, which differ in how they trade off route extension cost against the urgency of including each unrouted customer; the \textbf{savings} (Clarke--Wright) algorithm; a \textbf{nearest-neighbour} (NN) greedy strategy; and a \textbf{sweep} decoder.

\paragraph{HVASP decoder.} The HVASP decoder performs a two-phase assignment: first,
each activity is assigned to its eligible vessel with the best-matching random key,
subject to non-overlap and maintenance constraints; second, assigned activities are
sequenced within each vessel's schedule to respect time windows. The assignment is split into two phases to increase the likelihood that a chromosome decodes to a feasible schedule: each phase only needs to satisfy a subset of the constraints (eligibility, non-overlap, and maintenance in the first phase; time windows in the second), so the decoder resolves the combinatorially hardest decisions before committing to a temporal ordering, rather than requiring a single pass to satisfy all constraint families at once. Activities that cannot be feasibly assigned contribute to the normalised violation
$\bar{v}$ used by $p_\mathrm{cv}$.

\subsection{Instantiation II: IMAP-GA-II}
\label{subsec:imap-gaii}

IMAP-GA-II and IMAP-BRKGA are structurally opposite GA paradigms. BRKGA uses random-key chromosomes decoded by a problem-specific function, with feasibility handled through decoder-reported violation information and evolution proceeding through biased crossover of an explicit elite partition. IMAP-GA-II operates directly on real-valued decision vectors, with feasibility governed by constraint functions and evolution proceeding through standard real-coded genetic variation operators. Despite this structural opposition, both reduce to the same IMAP-IGS template (Algorithm~\ref{alg:imap-ga}): the current-best re-insertion (Section~\ref{subsec:intergenerational}) and the weighted constraint-violation preference function $p_\mathrm{cv}$ (Section~\ref{subsec:cv-pref}) are transferred identically to both instantiations, without modification to the variation operators. This host-agnostic transfer is the empirical demonstration that the IMAP fitness signal, not the choice of GA architecture, drives the framework's behaviour.

The name IMAP-GA-II is used deliberately to avoid presenting the method as an NSGA-II variant. NSGA-II is defined by its Pareto-front search machinery: fast non-dominated sorting and crowding-distance-based diversity preservation. Relative to the generic template of Algorithm~\ref{alg:imap-ga} and to standard NSGA-II, the differences are set in italics: \emph{fast non-dominated sorting and crowding-distance selection are removed entirely}. What remains is not Pareto-based selection, but a conventional real-coded GA operator stack: simulated binary crossover (SBX), polynomial mutation, population replacement, and \emph{tournament selection driven by a scalar fitness value} \cite{deb2002fast}, where \emph{the scalar fitness is the second-pass IMAP score computed within the current evaluation pool}.

Formally, IMAP-GA-II specialises IMAP-IGS with \textsc{EvolvePopulation} implemented by real-coded GA operators. Parent selection is driven solely by IMAP score, including the weighted constraint-violation preference when infeasible candidates are present. Offspring are generated using SBX crossover with parameter $\eta_c$ and polynomial mutation with parameter $\eta_m$. No non-dominated sorting, crowding distance, reference vectors, hypervolume indicators, or other Pareto-front preservation mechanisms are used. Consequently, IMAP-GA-II does not attempt to approximate a Pareto front to then select from it. It directly searches the continuous decision space for the feasible solution with maximum aggregated group preference. The current-best re-insertion from Section~\ref{subsec:intergenerational} applies identically, since the non-stationarity problem is a property of z-normalisation and is independent of the host GA's operators.

\section{Experimental Setup}
\label{sec:experiments}

This section defines the experimental protocol: the baseline algorithms and paired controls (Section~\ref{subsec:baselines}), the fairness rationale of the comparison (Section~\ref{subsec:fairness}), the conflicting and aligned stakeholder configurations that separate conflict resolution from generic optimisation performance (Section~\ref{subsec:stakeholder-configs}), the algorithm parameters and stopping criteria (Section~\ref{subsec:parameters}), and the evaluation metrics (Section~\ref{subsec:metrics}).

\subsection{Baseline algorithms}
\label{subsec:baselines}
Baselines fall into two categories: \emph{paired controls} that isolate the IMAP
signal from GA structure where such a paired control was implemented, and
\emph{state-of-the-art CMOEAs} representing the current standard for each problem
class.

\paragraph{Paired control.} BRKGA-WS uses the BRKGA structure of IMAP-BRKGA but
replaces the IMAP fitness signal with a weighted-sum of normalised objective values
using equal weights. Because BRKGA-WS shares the decoder, chromosome representation,
population structure, and evolutionary operators with IMAP-BRKGA, any performance
difference in the combinatorial domains is directly attributable to the choice between
IMAP and weighted-sum aggregation. No weighted-sum NSGA-II control is included for
DAS-CMOP; in that domain, the comparison is between IMAP-GA-II and standard
Pareto-based CMOEAs.

\paragraph{DAS-CMOP baselines.} Standard NSGA-II \cite{deb2002fast}, NSGA-III
\cite{deb2014evolutionary}, C-TAEA \cite{li2019ctaea}, and MOEA/D-CDP \cite{Fan2016moea/d-cdp}
are implemented using the pymoo library \cite{Blank2020Pymoo} with parameters
recommended in the respective publications.

\paragraph{MO-VRPTW baselines.} NSGA-III, SPEA2 \cite{zitzler2001spea2} and BRKGA-WS serve as Pareto-based routing baselines.

\paragraph{HVASP baselines.} Five CP solver variants establish the performance ceiling
of exact combinatorial methods for this problem class: CP-cost (single-objective,
minimise cost only), CP-duration (single-objective, minimise duration only), CP-hier
(hierarchical: minimise cost first, then minimise duration within 2\% cost slack),
CP-RWS (raw weighted-sum with equal weights over all model objectives), and CP-NWS
(normalised weighted-sum with equal weights). Each sub-solve is allocated a wall-clock
limit of 300\,s for Small and 600\,s for Medium instances. Pareto-based NSGA-II is additionally included as an evolutionary
baseline; its best representative is selected post-hoc via IMAP for the $Z^*$
comparison.

\paragraph{Post-hoc preference enrichment of baselines.} When a baseline returns a
set of candidate solutions, especially a Pareto-front approximation, the evaluation
does not force the baseline to use its native output directly. Instead, the stated PFM
preference model is applied post-hoc to the returned set and the best available
representative is selected for comparison. This deliberately enriches the baselines:
they receive preference information at the decision stage even though they did not use
it during search. The experiments therefore test whether using the preference model
during search improves the final decision beyond what can be recovered by applying the
same preference model only after search.

\subsection{Fairness of Comparison}
\label{subsec:fairness}
A possible concern is that IMAP-IGS is advantaged because it optimises the same quantity used for evaluation. This concern would apply if the experimental task were defined as generic Pareto-front approximation. However, the task considered in this paper is different: the decision model contains explicit stakeholder preference functions and the goal is to identify the feasible solution with highest group preference. Therefore, the post-hoc IMAP score is not an algorithm-specific metric but the decision criterion induced by the problem formulation itself. Pareto-based and weighted-sum methods are compared fairly because they are alternative search procedures for the same design--decision problem. Their final outputs are reduced to a single representative solution and evaluated according to the same preference model.

The experiment therefore compares search policies, not evaluation metrics. IMAP-IGS uses the stated preference model during search. Pareto-based methods first search for a non-dominated set and are then allowed to select the most preferred point from that set. Weighted-sum methods search for a scalarised objective compromise and are then evaluated against the same final decision criterion. In each case, the returned point is judged by the original question: which feasible point gives the highest group preference?

\begin{table}[t]
    \centering
    \caption{Comparison protocol: different search policies, same final decision criterion.}
    \label{tab:search_policy_comparison}
    \small
    \begin{tabularx}{\linewidth}{@{}p{0.22\linewidth}XX@{}}
        \toprule
        \textbf{Method family} & \textbf{Search policy} & \textbf{Final representative used for $Z^*$} \\
        \midrule
        IMAP-IGS & Directly searches for the feasible point with highest IMAP group preference. & Current-best IMAP solution $\mathbf{x}^\dagger$. \\
        Pareto-based CMOEAs & Approximate a non-dominated set without using the preference model during selection. & Most preferred feasible member of the returned set, selected post-hoc by IMAP. \\
        Weighted-sum baselines & Search for a scalar objective compromise using normalised objective values. & Best feasible representative returned by the scalarised run, evaluated by IMAP. \\
        CP baselines & Optimise cost, duration, hierarchy, or scalarised CP objectives. & Feasible CP solution returned by the solver, evaluated by IMAP. \\
        \bottomrule
    \end{tabularx}
\end{table}

\begin{figure}[t]
    \centering
    \includegraphics[width=0.95\linewidth]{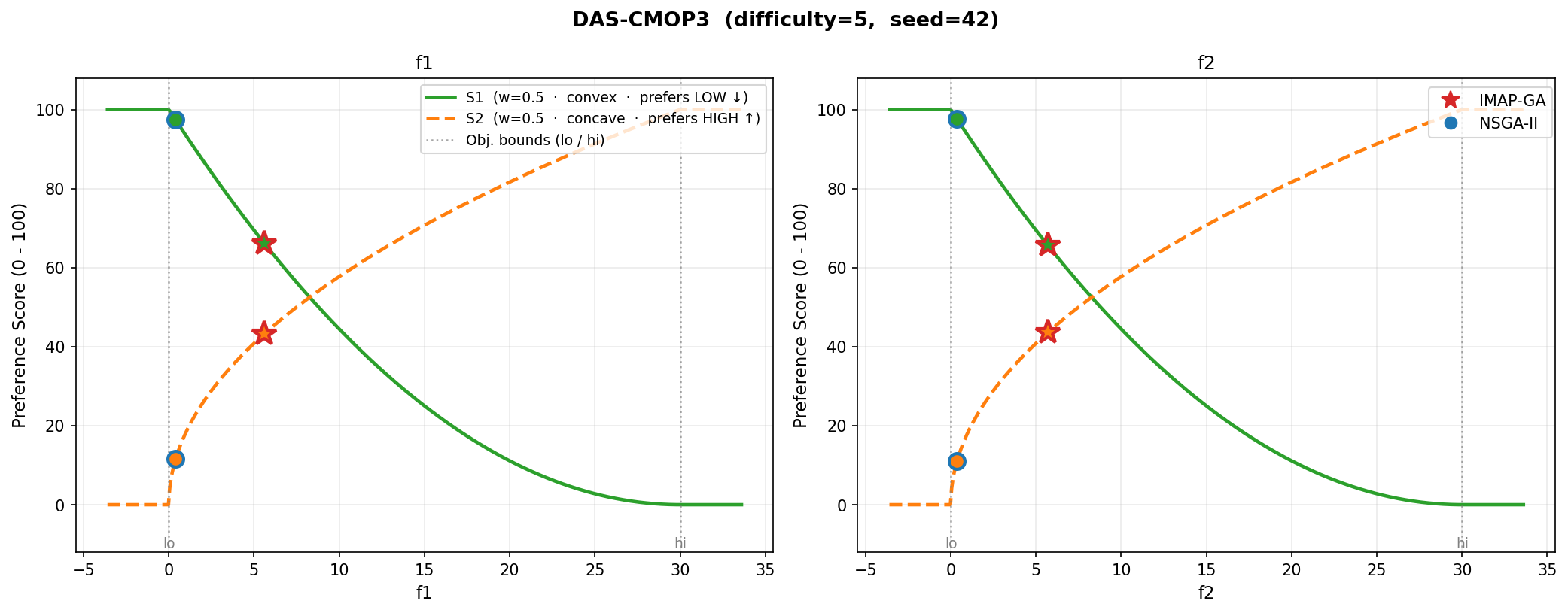}
    \caption{Illustrative DAS-CMOP3 conflicting-preference example with equal stakeholder importance weights ($W_1=W_2=0.5$). Stakeholder~1 prefers lower objective values, while Stakeholder~2 prefers higher objective values. IMAP-IGS returns a compromise point between the two preference directions, whereas the Pareto-based NSGA-II baseline remains near the raw-objective low region. The curves are illustrative preference functions used to visualise the trade-off mechanism; the benchmark configurations used for aggregate results are defined in Section~\ref{subsec:stakeholder-configs}. The example visualises why the relevant decision criterion is the point with highest aggregated preference, not the point favoured by Pareto-front optimisation alone.}
    \label{fig:dascmop3_preferences}
\end{figure}

Figure~\ref{fig:dascmop3_preferences} illustrates this distinction on a single DAS-CMOP3 run. For this illustrative plot, the two stakeholders are assigned equal importance weights ($W_1=W_2=0.5$) so that the desired decision is visually a compromise between opposing preference directions. The Pareto-based baseline moves toward the raw-objective low region because this is favoured by Pareto-front optimisation under the minimisation formulation. IMAP-IGS instead uses the preference functions during search and returns an interior trade-off point with substantially higher aggregate stakeholder preference. The figure is illustrative rather than a replacement for the full benchmark protocol; the aggregate results in Section~\ref{sec:results} use the stakeholder configurations defined in Section~\ref{subsec:stakeholder-configs}.

\subsection{Stakeholder Configurations}
\label{subsec:stakeholder-configs}

Every problem is evaluated under two configurations:

\textbf{Configuration C (Conflicting).} Stakeholder 1 has preference functions in the standard minimisation direction; Stakeholder 2 has reversed preference functions on the same objectives. Stakeholder importance weights are those defined for each problem family in Section~\ref{sec:problem_formulation}. This creates genuine inter-stakeholder conflict: no solution simultaneously satisfies both stakeholders at their respective optima, and neither Pareto-based nor weighted-sum methods natively possess a mechanism for targeting the specific group-preference compromise.

\textbf{Configuration A (Aligned).} Both stakeholders have preference functions in
the same minimisation direction with identical objective weights; importance weights
remain the same as in Configuration C. This eliminates preference conflict: the group-
optimal IMAP solution collapses onto the uniformly good region of the Pareto front,
which any method finding a good Pareto-efficient point can in principle identify.
Configuration A serves as a mechanism and fairness test: PFM theory predicts IMAP's advantage
should disappear when there is no conflict to resolve, because aligned monotone preferences make the preference-optimal point recoverable in principle by Pareto-front approximation followed by post-hoc IMAP selection (Appendix~\ref{appendix:pareto-imap}).

Problem-specific importance and objective weights are specified with each problem formulation. DAS-CMOP1--6 use stakeholder importance weights $(0.3,0.7)$ and equal objective weights $[0.5,0.5]$; DAS-CMOP7--9 and MO-VRPTW use stakeholder importance weights $(0.6,0.4)$, with $S_1$ using equal objective weights and $S_2$ placing greater emphasis on the third objective; HVASP uses stakeholder importance weights $(0.3,0.7)$ with equal objective weights $[0.5,0.5]$ over cost and duration. The illustrative plot in Figure~\ref{fig:dascmop3_preferences} deliberately uses equal stakeholder importance weights only to visualise the compromise mechanism; it is not the weighting protocol for the aggregate benchmark tables.

\subsection{Parameters and stopping criteria}
\label{subsec:parameters}

IMAP-BRKGA uses domain-tuned settings with a crossover bias of $\rho_b = 0.70$ in
both combinatorial domains. On HVASP it uses $|\mathcal{P}| = 300$, elite fraction
$\rho_e = 0.05$, and mutant fraction $\rho_m = 0.30$; experiments run for 150
generations or a 15-minute wall-clock cap per algorithm per instance, whichever is
reached first. On MO-VRPTW it uses $|\mathcal{P}| = 200$, $\rho_e = 0.15$, and
$\rho_m = 0.10$; every algorithm is stopped after 10{,}000 function evaluations or a
10-minute wall-clock cap, whichever comes first, so all MO-VRPTW baselines share the
same budget. DAS-CMOP experiments are stopped at 300{,}000 function evaluations,
matching the evaluation budgets of all baseline CMOEAs to ensure fair comparison.
IMAP-GA-II uses $|\mathcal{P}| = 300$ with SBX crossover
($\eta_c = 15$, rate 0.9) and polynomial mutation ($\eta_m = 20$, rate $1/n$).
Selection is based on IMAP score rather than Pareto rank; no non-dominated sorting,
crowding distance, or other Pareto-front maintenance mechanism is used.

\subsection{Performance metrics}
\label{subsec:metrics}
The primary metric is $Z^*$: the IMAP score of an algorithm's best feasible representative, computed post-hoc over the combined pool of all algorithms' best representatives. The winning algorithm's representative scores 100 by construction, all others are scored relative to it on a $[0, 100]$ scale. Separate $Z^*$ values are computed under Configuration C and Configuration A using the respective preference functions. Absolute scores are not comparable across configurations, but algorithm rankings within each configuration are the quantities of interest. The constraint violation preference function is excluded from the IMAP score calculation.

For Pareto-based algorithms, the returned non-dominated set is not discarded. Instead, the most preferred representative from that set is selected using the same preference model. Thus Pareto-based methods are not penalised for returning a set rather than a single point; they are given the benefit of post-hoc preference selection over their own returned approximation. The comparison therefore asks whether guiding the search by preferences produces a better final decision than applying the same preference model only after Pareto search has terminated.

Statistical reporting is based on win frequency. For DAS-CMOP, each (problem, difficulty, seed) triple is scored independently: the best feasible representative from each algorithm is pooled into a single IMAP evaluation call using theoretical (fixed) objective bounds. The algorithm with the highest $Z^*$ in that pool wins the seed. Win counts across 30 seeds are reported per (problem, difficulty) level, win frequency is the primary evidence of directional dominance. For MO-VRPTW, win counts are reported across three seeds per Solomon instance, aggregated per structural category (C1, C2, R1, R2, RC1, RC2). HVASP results are reported as win counts across the five instances per size class (Small and Medium); no formal hypothesis test is applied, as HVASP serves as structured real-world validation rather than a statistical benchmark.

\section{Experimental Results}
\label{sec:results}

\subsection{DAS-CMOP Results}
As introduced in Sections~\ref{subsec:das-cmop} and~\ref{sec:literature_review}, DAS-CMOP is used here as the continuous constraint-stress benchmark. The results in this subsection therefore focus on whether preference-guided search can identify the preferred feasible representative under increasing constraint difficulty, rather than re-describing the benchmark itself.

Table~\ref{tab:dascmop_meanwinrate_c} reports the mean win rate per problem, where win rate is the number of seeds (out of 30) in which each algorithm achieved the highest IMAP preference score. Values shown are averaged over all 16 difficulty levels. Table~\ref{tab:dascmop_wins_c} shows the per-seed win count across all (problem, difficulty) instances under Configuration C. Dashes indicate the algorithm was not run for that problem type. Because seeds in which no algorithm returns a feasible solution have no winner, and seeds in which several algorithms return equally preferred representatives are credited to each, per-problem win rates need not sum to 30. IMAP-GA-II wins the majority of seeds across all nine DAS-CMOP problems, narrowly on DAS-CMOP9. The advantage is largest at lower difficulty, where the feasible region is larger and IMAP-GA-II has fewer problems with infeasibility. At higher difficulties, the advantage becomes lower, as IMAP-GA-II struggles with the smaller feasible region. Table~\ref{tab:dascmop_meanwinrate_nc} reports the corresponding mean win rate under Configuration A (aligned preferences). The per-seed win count across all (problem, difficulty) instances under Configuration A is found in table~\ref{tab:dascmop_wins_a}.

Three findings emerge:

\textbf{Finding 1 (IMAP selection vs.\ Pareto selection).}
IMAP-GA-II uses the same standard real-coded variation operators often used in NSGA-II implementations, namely SBX crossover and polynomial mutation, but it removes the defining Pareto-front mechanisms: fast non-dominated sorting and crowding-distance selection. Standard NSGA-II consistently underperforms IMAP-GA-II under Configuration C. Pareto-based search explores the entire front, including regions irrelevant to the stakeholders' preferences under the reversed objective directions of Configuration C. Embedding preference maximisation from the first generation concentrates search on the conflict-resolving compromise, avoiding wasted evaluations in the wrong part of objective space. The comparison therefore attributes the performance gap to replacing Pareto-front selection with IMAP-guided selection, while keeping the underlying real-coded genetic variation operators standard.

\textbf{Finding 2 (Structural indirection, not NSGA-II weakness).}
C-TAEA and MOEA/D-CDP also underperform relative to IMAP-GA-II under Configuration C, despite being state-of-the-art CMOEAs designed specifically to handle constrained settings. This confirms that the relevant distinction is structural: methods built around Pareto-front approximation defer preference resolution to a post-hoc step, and this deferral is costly under genuine inter-stakeholder conflict regardless of how well the algorithm handles constraints.

\textbf{Finding 3 (Mechanism test: aligned preferences collapse the gap).}
Under Configuration A, the win-frequency advantage of IMAP-GA-II over standard NSGA-II, C-TAEA, and MOEA/D-CDP narrows substantially. The algorithms, instances, and evaluation budgets are unchanged between Configuration C and Configuration A; only the preference functions change. This collapse is the decisive mechanism test: IMAP's advantage is specifically attributable to its conflict-resolution mechanism, not to generic optimisation superiority. When there is no inter-stakeholder conflict to resolve, preference-guided and Pareto-based search converge on the same solutions.

\begin{longtable}{@{} l r r r r r @{}}
\caption{ Mean win rate per DAS-CMOP problem under Configuration C (conflicting preferences). }
\label{tab:dascmop_meanwinrate_c} \\
\toprule
\textbf{Problem} &
\textbf{IMAP-GA-II} &
\textbf{C-TAEA} &
\textbf{MOEA/D-CDP} &
\textbf{NSGA-II} &
\textbf{NSGA-III} \\
\midrule
\endfirsthead
\multicolumn{6}{c}{\tablename~\thetable{} \emph{(continued)}} \\[2pt]
\toprule
\textbf{Problem} &
\textbf{IMAP-GA-II} &
\textbf{C-TAEA} &
\textbf{MOEA/D-CDP} &
\textbf{NSGA-II} &
\textbf{NSGA-III} \\
\midrule
\endhead
\bottomrule
\endlastfoot
\begin{filecontents*}{data/das_cmop_conflict_mean_winrate_latex.csv}
Problem,IMAP_GA,C_TAEA,MOEA_D_CDP,NSGA_II,NSGA_III,_group_break
DAS-CMOP1,\textbf{21.6},0.0,0.9,1.9,--,0
DAS-CMOP2,\textbf{25.6},0.1,2.1,3.1,--,0
DAS-CMOP3,\textbf{24.6},0.1,1.6,3.2,--,0
DAS-CMOP4,\textbf{21.8},3.1,2.5,3.2,--,0
DAS-CMOP5,\textbf{21.1},1.7,5.2,2.8,--,0
DAS-CMOP6,\textbf{24.4},1.2,0.5,4.1,--,0
DAS-CMOP7,\textbf{24.4},0.3,1.4,--,3.9,1
DAS-CMOP8,\textbf{19.8},4.3,5.1,--,1.9,0
DAS-CMOP9,\textbf{15.5},0.2,11.8,--,3.6,0
\end{filecontents*}
\csvreader[
head,
late after line = \\,
]{data/das_cmop_conflict_mean_winrate_latex.csv}{}
{%
\ifnum\csvcolvii=1 \midrule \fi
\csvcoli & \csvcolii & \csvcoliii & \csvcoliv & \csvcolv & \csvcolvi
}
\end{longtable}

\begin{longtable}{@{} l r r r r r @{}}
\caption{Mean win rate per DAS-CMOP problem under Configuration A (aligned preferences).}
\label{tab:dascmop_meanwinrate_nc} \\
\toprule
\textbf{Problem} &
\textbf{IMAP-GA-II} &
\textbf{C-TAEA} &
\textbf{MOEA/D-CDP} &
\textbf{NSGA-II} &
\textbf{NSGA-III} \\
\midrule
\endfirsthead
\multicolumn{6}{c}{\tablename~\thetable{} \emph{(continued)}} \\[2pt]
\toprule
\textbf{Problem} &
\textbf{IMAP-GA-II} &
\textbf{C-TAEA} &
\textbf{MOEA/D-CDP} &
\textbf{NSGA-II} &
\textbf{NSGA-III} \\
\midrule
\endhead
\bottomrule
\endlastfoot
\begin{filecontents*}{data/das_cmop_no_conflict_mean_winrate_latex.csv}
Problem,IMAP_GA,C_TAEA,MOEA_D_CDP,NSGA_II,NSGA_III,_group_break
DAS-CMOP1,5.1,\textbf{8.4},4.1,7.2,--,0
DAS-CMOP2,\textbf{13.1},2.9,4.4,10.4,--,0
DAS-CMOP3,\textbf{13.2},3.6,4.2,9.6,--,0
DAS-CMOP4,9.9,3.3,\textbf{14.4},3.9,--,0
DAS-CMOP5,5.1,4.6,6.3,\textbf{14.5},--,0
DAS-CMOP6,\textbf{13.2},4.2,11.0,2.3,--,0
DAS-CMOP7,\textbf{11.1},3.4,7.6,--,7.9,1
DAS-CMOP8,\textbf{15.2},3.0,6.5,--,5.4,0
DAS-CMOP9,\textbf{13.6},7.8,6.8,--,2.3,0
\end{filecontents*}
\csvreader[
head,
late after line = \\,
]{data/das_cmop_no_conflict_mean_winrate_latex.csv}{}
{%
\ifnum\csvcolvii=1 \midrule \fi
\csvcoli & \csvcolii & \csvcoliii & \csvcoliv & \csvcolv & \csvcolvi
}
\end{longtable}

\subsection{MO-VRPTW Results}
The Solomon benchmark tests preference-guided routing across six structural categories differing in time-window tightness, clustering, and demand distribution. Under conflicting preferences over route distance, vehicle count, and route balance, the challenge is not only to find feasible routes but to identify the specific feasible route plan that best satisfies the group preference. Pareto-front methods are indirect for this task: the non-dominated set of routes may be large and distributed across objective space, with none of its members guaranteed to correspond to the preference-optimal compromise without additional selection logic that the algorithm was not designed to provide.

Table~\ref{tab:vrptw_c_mean} reports win rates across three seeds per Solomon instance, aggregated per structural category, under Configuration C. For each problem instance, win rate is the number of seeds (out of 3) in which each algorithm achieved the highest preference score. Values shown are averages over all instances within each category. Table~\ref{tab:vrptw_a_mean} reports the corresponding win under Configuration A. Table~\ref{tab:vrptw_conflict_decoder_ranking} shows that the savings heuristic performs best as decoder in all three tight-time-window categories (C100, R100, RC100), while no single decoder dominates the wide-window 200-series, where I1, sweep, and nearest-neighbour lead on C200, R200, and RC200 respectively. Across all decoder choices, the IMAP fitness signal consistently dominates the non-preference-guided alternatives.

Four findings emerge:

\textbf{Finding 1 (IMAP-BRKGA identifies the preference-optimal route plan).}
IMAP-BRKGA wins the majority of instances across all six Solomon categories under Configuration C. The advantage is strongest in R2 and RC2 categories, where wider time windows create a larger feasible space and the IMAP fitness signal has more room to distinguish between routes with different preference profiles.

\textbf{Finding 2 (Weighted-sum scalarisation misaligns under stakeholder conflict).}
IMAP-BRKGA outperforms BRKGA-WS in all six categories. The same decoder, crossover, and mutation produce a worse outcome when the fitness signal is a weighted-sum of normalised objectives rather than the IMAP score. Equal-weight weighted-sum scalarisation cannot target the asymmetric preference optimum that arises under Configuration C's reversed stakeholder directions and unequal importance weights.

\textbf{Finding 3 (Pareto representatives miss the preference optimum).}
The best IMAP-elected representatives from NSGA-II, NSGA-III and SPEA2 consistently score below IMAP-BRKGA's $\mathbf{x}^\dagger$. Pareto-based search produces non-dominated route plans that are efficient in objective space but are not guided toward the preference-optimal compromise. The post-hoc IMAP selection step, applied to a Pareto front that was not guided by preference, often does not recover the solution that preference-guided search finds, particularly when the preference optimum lies at a point in objective space that is far from the central tendency of the front.

\textbf{Finding 4 (Mechanism test: aligned preferences collapse the gap).}
Under Configuration A, the dominance of IMAP-BRKGA is less pronounced across almost all six categories. The same algorithms, the same instances, the same decoder: only the preference functions differ. The collapse of the performance gap confirms that IMAP's advantage on MO-VRPTW is specifically attributable to its conflict-resolution mechanism.

\begin{table}[h]
    \centering
    \footnotesize \caption{Mean win rate per Solomon instance category under Configuration C (conflicting preferences).}
    \label{tab:vrptw_c_mean}
    \footnotesize
    \begin{tabular}{@{} l r r r r r @{}}
        \toprule
        \textbf{Category} &
        \textbf{IMAP-BRKGA} &
        \textbf{BRKGA-WS} &
        \textbf{NSGA-II} &
        \textbf{NSGA-III} &
        \textbf{SPEA2} \\
        \midrule
        \begin{filecontents*}{data/vrptw_conflict_mean_winrate_latex.csv}
Category,IMAP_BRKGA,BRKGA_WS,NSGA_II,NSGA_III,SPEA2,_group_break
C100,\textbf{2.9},0.0,0.1,0.0,0.0,0
C200,\textbf{2.9},0.0,0.5,0.4,0.4,0
R100,\textbf{2.9},0.0,0.0,0.0,0.1,1
R200,\textbf{3.0},0.0,0.0,0.0,0.0,0
RC100,\textbf{2.6},0.1,0.1,0.0,0.1,1
RC200,\textbf{3.0},0.0,0.0,0.0,0.0,0
\end{filecontents*}
\csvreader[
        head,
        late after line = \\,
        ]{data/vrptw_conflict_mean_winrate_latex.csv}{}
        {%
        \ifnum\csvcolvii=1 \midrule \fi
        \csvcoli & \csvcolii & \csvcoliii & \csvcoliv & \csvcolv & \csvcolvi
        }
        \bottomrule
    \end{tabular}
\end{table}

\begin{table}[h]
    \centering
    \footnotesize \caption{Mean win rate per Solomon instance category under Configuration A (aligned preferences).}
    \label{tab:vrptw_a_mean}
    \footnotesize
    \begin{tabular}{@{} l r r r r r @{}}
        \toprule
        \textbf{Category} &
        \textbf{IMAP-BRKGA} &
        \textbf{BRKGA-WS} &
        \textbf{NSGA-II} &
        \textbf{NSGA-III} &
        \textbf{SPEA2} \\
        \midrule
        \begin{filecontents*}{data/vrptw_no_conflict_mean_winrate_latex.csv}
Category,IMAP_BRKGA,BRKGA_WS,NSGA_II,NSGA_III,SPEA2,_group_break
C100,\textbf{2.7},0.2,0.2,0.2,0.0,0
C200,\textbf{2.0},\textbf{2.0},0.0,0.0,0.2,0
R100,\textbf{1.8},1.0,0.7,0.3,0.5,1
R200,\textbf{2.7},0.4,0.0,0.2,0.0,0
RC100,1.2,\textbf{2.2},0.2,0.2,0.8,1
RC200,\textbf{2.8},0.0,0.5,0.0,0.0,0
\end{filecontents*}
\csvreader[
        head,
        late after line = \\,
        ]{data/vrptw_no_conflict_mean_winrate_latex.csv}{}
        {%
        \ifnum\csvcolvii=1 \midrule \fi
        \csvcoli & \csvcolii & \csvcoliii & \csvcoliv & \csvcolv & \csvcolvi
        }
        \bottomrule
    \end{tabular}
\end{table}

\begin{table}[ht]
  \centering
  \footnotesize  \caption{Mean IMAP score per BRKGA-IMAP decoder per instance category on the MO-VRPTW conflicting benchmark, averaged across all seeds.}
  \label{tab:vrptw_conflict_decoder_ranking}
  \footnotesize
  \begin{tabular}{@{} l r r r r r r r @{}}
      \toprule
      \textbf{Category} &
      \textbf{I1} &
      \textbf{I2} &
      \textbf{I3} &
      \textbf{NN} &
      \textbf{Savings} &
      \textbf{Sweep} &
      \textbf{Random} \\
      \midrule
      \begin{filecontents*}{data/vrptw_conflict_decoder_ranking_latex.csv}
Category,(i1),(i2),(i3),(nn),(random),(savings),(sweep),_group_break
C100,78.9,70.6,60.2,91.6,85.3,\textbf{98.4},3.3,0
C200,\textbf{94.7},62.4,65.3,80.7,73.1,75.8,4.7,0
R100,78.7,43.2,8.8,68.1,74.8,\textbf{85.5},65.6,1
R200,65.5,3.0,49.1,77.3,61.9,72.3,\textbf{93.9},0
RC100,67.1,58.6,3.1,66.3,66.6,\textbf{80.0},46.0,1
RC200,62.8,14.6,38.7,\textbf{89.3},73.9,83.5,80.1,0
\end{filecontents*}
\csvreader[
      head,
      late after line = \\,
      ]{data/vrptw_conflict_decoder_ranking_latex.csv}{}
      {
      \ifnum\csvcolix=1 \midrule \fi
      \csvcoli & \csvcolii & \csvcoliii & \csvcoliv & \csvcolv & \csvcolvi & \csvcolvii & \csvcolviii
      }
      \bottomrule
  \end{tabular}
\end{table}

\subsection{HVASP results}
Table~\ref{tab:hvasp_results_c} reports the win count per size class for all eight algorithms grouped by size class under Configuration C. Table~\ref{tab:hvasp_results_a} reports the corresponding results under Configuration A. Four findings emerge:

\textbf{Finding 1 (IMAP-BRKGA scales to real-world scheduling).}
IMAP-BRKGA achieves the highest $Z^*$ in both size classes and in all individual instances under Configuration C. The advantage is maintained from Small to Medium instances, indicating that preference-guided search remains effective as vessel counts, activity counts, and constraint complexity grow.

\textbf{Finding 2 (IMAP signal over BRKGA structure).}
IMAP-BRKGA outperforms BRKGA-WS across all size classes with the same decoder, the same chromosome representation, and the same operators. The difference is the fitness signal: IMAP aggregation identifies the preference-optimal assignment, weighted-sum aggregation does not reliably do so under conflicting preference directions.

\textbf{Finding 3 (Pareto selection insufficient for combinatorial scheduling).}
NSGA-II's best IMAP-selected representative consistently scores below IMAP-BRKGA's $\mathbf{x}^\dagger$. The combinatorial feasibility structure of HVASP is not well served by NSGA-II's real-valued operators, and Pareto selection does not guide search toward the specific cost-makespan compromise that the preference model requires.

\textbf{Finding 4 (Mechanism test: aligned preferences collapse the gap).} Under Configuration A, IMAP-BRKGA and BRKGA-WS converge in $Z^*$ across all size classes. NSGA-II's IMAP-selected representative becomes competitive and CP-RWS and CP-NWS close the gap substantially. The collapse of IMAP-BRKGA's advantage confirms that its superiority under Configuration C is attributable to conflict-resolution, not to incidental optimisation superiority over the combinatorial search problem itself.

\begin{table}[h]
    \centering
    \footnotesize \caption{\small Win count per size class under Configuration C (conflicting preferences).}
    \label{tab:hvasp_results_c}
    \footnotesize
    \begin{tabular}{@{} l r r r r r r r r @{}}
    \toprule
    \textbf{Size} &
    \textbf{IMAP-BRKGA} &
    \textbf{BRKGA-WS} &
    \textbf{NSGA-II} &
    \textbf{CP-cost} &
    \textbf{CP-dur} &
    \textbf{CP-Hier} &
    \textbf{CP-RWS} &
    \textbf{CP-NWS} \\
    \midrule
    \begin{filecontents*}{data/hvasp_conflict_winrate_latex.csv}
Size,BRKGA_IMAP,BRKGA_WS,NSGA_II,CP_cost,CP_duration,CP_Hier,CP_RWS,CP_NWS,_group_break
Small,\textbf{5.0},0.0,0.0,0.0,0.0,0.0,0.0,0.0,0
Medium,\textbf{5.0},0.0,0.0,0.0,0.0,0.0,0.0,0.0,1
\end{filecontents*}
\csvreader[
    head,
    late after line = \\,
    ]{data/hvasp_conflict_winrate_latex.csv}{}
    {
    \csvcoli & \csvcolii & \csvcoliii & \csvcoliv & \csvcolv &
    \csvcolvi & \csvcolvii & \csvcolviii & \csvcolix
    }
    \bottomrule
    \end{tabular}
\end{table}

\begin{table}[h]
    \centering
    \footnotesize \caption{\small Win count per size class under Configuration A (aligned preferences).}
    \label{tab:hvasp_results_a}
    \footnotesize
    \begin{tabular}{@{} l r r r r r r r r @{}}
    \toprule
    \textbf{Size} &
    \textbf{IMAP-BRKGA} &
    \textbf{BRKGA-WS} &
    \textbf{NSGA-II} &
    \textbf{CP-cost} &
    \textbf{CP-dur} &
    \textbf{CP-Hier} &
    \textbf{CP-RWS} &
    \textbf{CP-NWS} \\
    \midrule
    \begin{filecontents*}{data/hvasp_no_conflict_winrate_latex.csv}
Size,BRKGA_IMAP,BRKGA_WS,NSGA_II,CP_cost,CP_duration,CP_Hier,CP_RWS,CP_NWS,_group_break
Small,\textbf{5.0},\textbf{5.0},2.0,4.0,0.0,0.0,\textbf{5.0},\textbf{5.0},0
Medium,3.0,3.0,0.0,3.0,0.0,0.0,\textbf{4.0},1.0,1
\end{filecontents*}
\csvreader[
    head,
    late after line = \\,
    ]{data/hvasp_no_conflict_winrate_latex.csv}{}
    {
    \csvcoli & \csvcolii & \csvcoliii & \csvcoliv & \csvcolv &
    \csvcolvi & \csvcolvii & \csvcolviii & \csvcolix
    }
    \bottomrule
    \end{tabular}
\end{table}

\subsection{Scaling to Many Objectives}
\label{subsec:results-dascmaop}
A structural property of IMAP-IGS is that the number of objectives, and hence of stakeholder preference functions, is immaterial to the optimisation process: the affine aggregation collapses any number of z-normalised preference scores into a single scalar fitness, so the framework requires no modification, no reference directions, no decomposition weights, and no archive management as the objective count grows. Pareto-based selection, by contrast, loses selective pressure in the many-objective regime, where almost all solutions become mutually non-dominated. To substantiate this claim empirically, the identical, unmodified IMAP-GA-II was additionally run on DAS-CMaOP1--9, the many-objective counterparts of DAS-CMOP defined by the same toolkit \cite{Fan_2020_DASCMOP}, with $m \in \{5, 8, 10\}$ objectives, over twelve difficulty settings and ten seeds per (problem, difficulty, $m$) cell, against the many-objective-capable baselines NSGA-III, C-TAEA, and MOEA/D-CDP at identical per-$m$ evaluation budgets.

Table~\ref{tab:dascmaop_summary} summarises the outcome. Under Configuration C, IMAP-GA-II attains the highest win rate in 26 of the 27 (problem, $m$) cells, and its mean win share rises with the objective count while every Pareto-based baseline collapses to single digits at $m \geq 8$. Under Configuration A the gap narrows but, unlike in the bi- and tri-objective results above, does not close: at $m \geq 5$ the returned front approximations are sparse and degenerate, so post-hoc selection frequently has no point near the preference optimum to choose. The purpose of this summary is solely to demonstrate that the method carries over to many objectives without any algorithmic change. The complete per-problem win-rate tables are provided in Appendix~\ref{appendix:dascmaop}; the accompanying feasibility, computational-cost, and per-stakeholder analyses are reported in the thesis version of this work.

\begin{table}[h]
    \centering
    \footnotesize \caption{Mean win share (\%) per algorithm and objective count $m$ on DAS-CMaOP, averaged over the nine problems, under Configuration C (conflicting preferences) and Configuration A (aligned preferences).}
    \label{tab:dascmaop_summary}
    \footnotesize
    \begin{tabular}{@{} l r r r r r r @{}}
        \toprule
        & \multicolumn{3}{c}{\textbf{Configuration C}} & \multicolumn{3}{c}{\textbf{Configuration A}} \\
        \cmidrule(lr){2-4} \cmidrule(l){5-7}
        \textbf{Algorithm} & $\mathbf{m{=}5}$ & $\mathbf{m{=}8}$ & $\mathbf{m{=}10}$ & $\mathbf{m{=}5}$ & $\mathbf{m{=}8}$ & $\mathbf{m{=}10}$ \\
        \midrule
        \begin{filecontents*}{data/dascmaop_mean_winshare_latex.csv}
Algorithm,C_m5,C_m8,C_m10,A_m5,A_m8,A_m10
IMAP-GA-II,67.0,77.7,74.8,48.4,40.4,43.0
NSGA-III,5.4,1.5,2.5,20.3,21.4,21.7
C-TAEA,15.3,5.6,6.7,10.0,8.1,6.0
MOEA/D-CDP,4.6,0.5,0.4,13.7,15.3,14.7
\end{filecontents*}
\csvreader[
        head,
        late after line = \\,
        ]{data/dascmaop_mean_winshare_latex.csv}{}
        {
        \csvcoli & \csvcolii & \csvcoliii & \csvcoliv & \csvcolv & \csvcolvi & \csvcolvii
        }
        \bottomrule
    \end{tabular}
\end{table}

\subsection{Discussion}
Three structurally distinct problem classes, a continuous constraint-stress-tested benchmark suite (DAS-CMOP), a combinatorial NP-hard routing problem (MO-VRPTW), and a real-world heterogeneous scheduling problem (HVASP), yield the same directional result: preference-guided search more reliably finds the preference-optimal solution than methods that use the preference model only after search under stakeholder conflict. Consistent success across this structural diversity rules out problem-specific explanations. The feasibility landscapes differ fundamentally, smooth continuous constraints, combinatorial time-window and capacity constraints, discrete assignment- and-scheduling constraints, yet the IMAP fitness signal succeeds across all three.

The IMAP fitness signal is the performance driver, established by two complementary comparisons. For combinatorial problems (MO-VRPTW, HVASP), BRKGA-WS shares every structural component with IMAP-BRKGA except the fitness signal. IMAP-BRKGA consistently outranks it. For DAS-CMOP, standard NSGA-II provides the Pareto-front baseline: it uses non-dominated sorting and crowding distance to approximate the Pareto front, after which the returned solutions are evaluated under the same preference model. IMAP-GA-II, by contrast, removes these Pareto-specific selection mechanisms and uses the same class of real-coded genetic variation operators, such as SBX crossover and polynomial mutation, with IMAP preference scoring as the selection signal. The fact that IMAP-guided search outperforms both the BRKGA weighted-sum control in the combinatorial domains and the Pareto-front baseline in the continuous DAS-CMOP domain supports the same conclusion across structurally different GA architectures: the observed advantage is driven by preference-guided fitness, not by a host-specific implementation detail.

The aligned-preference control (Configuration A) is the decisive verification and an empirical fairness check. Under aligned monotone preferences, Pareto-based methods are expected to become competitive because the preference-optimal solution lies in the same efficient region they are designed to approximate, provided that region is represented in the returned front. Baselines that were clearly dominated under Configuration C become competitive under Configuration A without any change to algorithms, instances, or evaluation protocol. This indicates that the advantage of IMAP-IGS is not generic optimisation superiority but stems from its conflict-resolution mechanism -- the behaviour PFM theory predicts for direct preference aggregation. The empirical findings are therefore not merely a demonstration that IMAP-IGS performs well, but also a test of the underlying theoretical prediction, and the prediction holds across three structurally distinct problem classes.

\section{Discussion and Future Work}
\label{sec:discussion}

\subsection{Scope and limitations}

The results support the paper's central scope claim: IMAP-IGS is most useful when tight feasibility constraints and genuinely conflicting stakeholder preferences must be handled within one evolutionary search process. Pareto-front approximation remains useful when the purpose is exploration, ordering, or elicitation from a trade-off surface. IMAP-IGS addresses a different design--decision setting, where stakeholder preference functions and weights are part of the decision model and the required output is a single implementable solution with maximum aggregated group preference.

The experiments also show that the benefit is conditional. Under aligned preferences, the advantage of direct IMAP-guided search narrows because the preference-optimal point is easier to recover through a good Pareto-front approximation followed by IMAP selection, or through a scalar method whose search direction happens to match the shared preference direction. This is not a weakness of the framework; it is the expected mechanism check. The method is designed for the harder regime in which feasibility constraints and stakeholder preference conflict are coupled.

The constraint-violation preference function should likewise be interpreted as a search-time device. It contributes to the method by keeping feasibility steering inside the same preference-based fitness structure, especially when early populations are fully infeasible. It is not a substantive stakeholder objective, and it is excluded from the final reported $Z^*$ scores. This distinction is important: the final decision remains based on stakeholder preference functions over feasible representatives.

Current-best re-insertion is necessary because IMAP scores are z-normalised within the current evaluation pool. Pool-relative scoring is mathematically appropriate for the affine preference aggregation used here, but it creates intergenerational dynamics that ordinary elitism cannot handle by simply carrying over a numerical score. Current-best re-insertion preserves the best currently known feasible solution by objective values and re-evaluates it inside the live pool, making the comparison valid within the current normalisation context.

The framework depends on the supplied preference model. If stakeholder functions or weights are misspecified, IMAP-IGS will optimise the wrong decision model. This limitation is shared by preference-based decision-support methods. The relevant future question is therefore not whether IMAP-IGS can infer preferences without elicitation, but how robust the selected decision is to plausible perturbations of the preference model, its slopes, breakpoints, target values, and stakeholder importance weights: whether the returned decision changes under such perturbations, and how much aggregated preference is lost when it does.

\subsection{Future work}

\begin{itemize}
    \item \textbf{Preference robustness and elicitation}: perturb the preference model, its slopes, breakpoints, target values, and stakeholder importance weights, and quantify whether the returned decision $\mathbf{x}^\dagger$ changes and how much aggregated preference is lost when it does; integrate stakeholder-in-the-loop elicitation so that functions and weights can be refined during the IMAP loop.
    \item \textbf{Constraint relaxation for over-constrained instances}: use PFM to guide which constraints to relax, and to what extent, when no feasible solution exists.
    \item \textbf{Further GA instantiations}: test evolutionary strategies, differential evolution, and memetic algorithms to further establish the host-agnostic claim and identify which GA structural properties amplify IMAP's preference-alignment benefit.
    \item \textbf{Convergence-threshold sensitivity}: evaluate the two-pass ranking threshold $\tau$ across values such as $20$, $40$, and $60$ to quantify how strongly convergence speed and final $Z^*$ depend on the filtering rule.
    \item \textbf{Feasibility-success reporting}: report the fraction of runs in which each algorithm finds at least one feasible representative, especially on the highest-difficulty DAS-CMOP settings, to separate feasibility discovery from preference quality conditional on feasibility.
    \item \textbf{Scalability analysis}: extend the MO-VRPTW experiments to Homberger and Gehring instances with up to 1000 customers, and extend HVASP to larger fleets and activity portfolios.
    \item \textbf{Dynamic and stochastic extensions}: use partial re-optimisation from $\mathbf{x}^\dagger$ when new activities, disruptions, or updated vessel information arrive.
\end{itemize}

\section{Conclusion}
\label{sec:conclusion}

This paper presented IMAP-IGS, a metaheuristic framework for highly constrained multi-objective design--decision problems with stakeholder preference conflict. The framework embeds IMAP preference aggregation, grounded in preference function modelling, directly into the evolutionary fitness signal. Its central algorithmic property is that IMAP replaces the scalar fitness function, not the host algorithm's representation or variation operators. This makes the approach host-agnostic: the same preference-guided signal can be embedded in structurally different evolutionary architectures.

Two instantiations were demonstrated. IMAP-BRKGA targets combinatorial problems with hard feasibility structures, while IMAP-GA-II targets continuous constrained benchmark problems using a real-coded GA operator stack. IMAP-GA-II deliberately removes the Pareto-front components associated with NSGA-II, namely non-dominated sorting and crowding-distance selection, retaining standard genetic variation operators such as SBX crossover and polynomial mutation. Both instantiations use the same constraint-violation preference function for search-time feasibility steering and the same current-best re-insertion mechanism for stabilising pool-relative $z$-normalised scoring.

The experiments evaluated three structurally distinct problem classes under two stakeholder configurations: conflicting preferences (Configuration C) and aligned preferences (Configuration A). Under conflicting preferences, IMAP-BRKGA outperformed its paired weighted-sum BRKGA control in the combinatorial domains, while IMAP-GA-II outperformed Pareto-based CMOEA baselines on DAS-CMOP; a supplementary many-objective campaign on DAS-CMaOP showed that this advantage widens further as the objective count grows, with no modification to the framework. Under aligned preferences, the advantage narrowed substantially. This collapse is central to the interpretation of the results: it shows that the performance gain is not generic optimisation superiority, but the expected benefit of using a preference-based conflict-resolution mechanism when stakeholders genuinely disagree.

The constraint-violation preference function is a separate methodological contribution for highly constrained settings. It allows feasibility proximity to influence search through the same preference-based fitness structure, particularly when the live population is fully infeasible. At the same time, it remains a pseudo-preference used only during search. Final reported $Z^*$ scores are computed over feasible representatives and exclude the constraint-violation term, so the final decision remains based on substantive stakeholder preferences.

The comparison is therefore a comparison of search policies for the same design--decision problem. Pareto-based and scalarised methods remain useful in settings for which exploration, ordering, or simple scalar trade-offs are the intended output. IMAP-IGS is intended for the regime where the decision model specifies stakeholder preferences and the required output is a single implementable solution. The results suggest that preference-based evolutionary search is most valuable where real-world constrained decision problems are hardest: when tight feasibility constraints and conflicting stakeholder preferences must be resolved together.

\newpage
\printbibliography

\newpage
\appendix

\section{Heterogeneous Vessel Allocation and Scheduling Problem}
\label{appendix:HVASP}

The HVASP is defined over four index sets:

\begin{table}[h]
    \centering
    \begin{tabular}{l|l}
        \toprule
        \textbf{Set} & \textbf{Description} \\
        \midrule
        $V = \{v_0, v_1, \ldots, v_n\}$ & Set of $n$ vessels \\
        $A = \{a_0, a_1, \ldots, a_m \}$ & Set of $m$ activities, partitioned into towing activities $A^{\text{tow}} \subseteq A$\\
        & and maintenance activities $A^{\text{maint}} \subseteq A$\\
        $R = \{r_0, r_1, \ldots, r_p\}$ & Set of $p$ roles \\
        $L = \{\ell_0, \ell_1, \ldots, \ell_q\} $ & Set of $q$ locations \\
        \bottomrule
    \end{tabular}
    \caption{All sets in the fleet allocation problem.}
    \label{tab:sets}
\end{table}

\textbf{Decision vector.} The controllable decision vector $\mathbf{x} = (x_1, \ldots, x_6)$ has components:

\begin{table}[h!]
    \centering
    \small
    \caption{Domain constraints for the simplified decision vector $x = (x_1, \dots, x_6)$}
    \begin{tabular}{c l l}
    \toprule
    \textbf{$\mathbf{x}$} & \textbf{Description} & \textbf{$g_f^{(0)}(x_i)$} \\
    \midrule
    $x_1$ & Start time of activity $a \in A$ & $\underline{T}_a \le x_1[a] \le \overline{T}_a$ \\
    $x_2$ & Location choice for maintenance activity $a \in A_{maint}$ & $x_2[a] \in L_a$ \\
    $x_3$ & Vessel assigned to role $r \in R$ & $x_3[r] \in \mathcal{D}_r$ \\
    $x_4$ & Sequencing variable: role $r'$ follows $r$ & $x_4[r,r'] \in \{0,1\}$ \\
    $x_5$ & Sequence start indicator: role $r$ is first in sequence & $x_5[r] \in \{0,1\}$ \\
    $x_6$ & Average speed for sub-route $(r, r')$ & $x_6[r, r'] \in [\underline{s}, \overline{s}]$\\
    \bottomrule
    \end{tabular}
    \label{tab:domain_constraints_simplified}
\end{table}

\textbf{Parameter vector.} The exogenous parameter vector $\mathbf{y} = (y_1, \ldots, y_{15})$ is defined in Table \ref{tab:parameter_vector_constraints}.

\begin{table}[h]
    \centering
    \small
    \caption{Parameter vector $\mathbf{y} = (y_1, \dots, y_{15})$ and their domain constraints}
    \begin{tabularx}{\textwidth}{c l X}
    \toprule
    \textbf{$\mathbf{y}$} & \textbf{Description} & \textbf{$g_f^{(0)}(y_i)$} \\
    \midrule
    $y_1$ & Duration of activity $a$ & \(y_1[a] \ge 0\) (Days) \\
    $y_2$ & Start time window for activity $a$ & \(y_2[a] = [\underline{y_2[a]}, \overline{y_2[a]}] \ge 0\) \\
    $y_3$ & Start location for towing activity $a \in A_{tow}$ & \(y_3[a] \in \text{Locations}\) \\
    $y_4$ & End location for towing activity $a \in A_{tow}$ & \(y_4[a] \in \text{Locations}\) \\
    $y_5$ & Locations for maintenance activity $a \in A_{maint}$ & \(y_5[a] \subseteq \text{Locations}\) \\
    $y_6$ & Predecessor of activity $a$ & \(y_6[a] \in A \cup \{\emptyset\}\) \\
    $y_7$ & Parent activity of role $r$ & \(y_7[r] \in A\) \\
    $y_8$ & Vessel domain for role $r$ & \(y_8[r] \subseteq V\) \\
    $y_9$ & Set of roles belonging to activity $a$ & \(y_9[a] = \{r \in R \mid \alpha(r)=a\}\) \\
    $y_{10}$ & Sailing distance from location $\ell$ to $\ell'$ & \(y_{10}[\ell,\ell'] \ge 0\) \\
    $y_{11}$ & Daily mobilisation rate for vessel $v$ & \(y_{11}[v] \ge 0\) \\
    $y_{12}$ & Fuel consumption rate function for vessel $v$ &
    \parbox[t]{0.75\textwidth}{%
    \(\displaystyle y_{12}[v](s) \ge 0,\;\forall s \in [y_{13}[v]^{\min}, y_{13}[v]^{\max}]\)
    } \\
    $y_{13}$ & Feasible sailing speed range for vessel $v$ &
    \parbox[t]{0.75\textwidth}{%
    \(\displaystyle y_{13}[v] = (s_v^{(0)}, s_v^{(1)}, \dots, s_v^{(|S_v|-1)}),\; s_v^{(k)} > 0\)
    } \\
    $y_{14}$ & Fuel price for vessel $v$ & \(y_{14}[v] \ge 0\) \\
    $y_{15}$ & Standby cost discount factor for vessel $v$ & \(y_{15}[v] \in [0,1]\) \\
    \bottomrule
    \end{tabularx}
    \label{tab:parameter_vector_constraints}
\end{table}

The role start and end locations are determined by activity type:

\[
\ell^{\text{start}}(r) = \begin{cases} y_3[y_7[r]] & \text{if } y_7[r] \in \mathcal{A}^{\text{tow}} \\ x_2[y_7[r]] & \text{if } y_7[r] \in \mathcal{A}^{\text{maint}} \end{cases}, \qquad \ell^{\text{end}}(r) = \begin{cases} y_4[y_7[r]] & \text{if } y_7[r] \in \mathcal{A}^{\text{tow}} \\ x_2[y_7[r]] & \text{if } y_7[r] \in \mathcal{A}^{\text{maint}} \end{cases}.
\]

The sailing duration for vessel $v$ to travel from $\ell$ to $\ell'$ at average speed $s$ is:

\[
\theta(v, \ell, \ell', s) = \left\lceil \frac{y_{10}[\ell, \ell']}{24 \cdot s} \right\rceil
\]

The total transition time between role $r$ and $r'$ is $\delta_{r,r'} = x_1[y_7[r']] - (x_1[y_7[r]] + y_1[y_7[r]])$, and the cost for vessel $v$ to travel between consecutive roles is:

\[
\gamma(v, r, r', s, \delta_{r,r'}) = \theta\!\left(v, \ell^{\text{end}}(r), \ell^{\text{start}}(r'), s\right) \cdot \Bigl(y_{11}[v] + y_{14}[v] \cdot y_{12}[v](s) - y_{15}[v] \cdot y_{11}[v]\Bigr) + y_{15}[v] \cdot y_{11}[v] \cdot \delta_{r,r'}.
\]

\textbf{Objective functions.} The two objective functions are:

$$f_1(\mathbf{x}, \mathbf{y}) = \sum_{r \in \mathcal{R}} \sum_{r' \in \mathcal{R}} x_4[r,r'] \cdot \gamma\!\left(x_3[r],\ r,\ r',\ x_6[r,r'],\ \delta_{r,r'}\right) \quad \text{(total cost)},$$

$$f_2(\mathbf{x}, \mathbf{y}) = \max_{a \in \mathcal{A}}\!\left(x_1[a] + y_1[a]\right) - \min_{a \in \mathcal{A}} x_1[a] \quad \text{(make-span)}.$$

\textbf{Feasibility constraints.} The 13 constraints partition into activity, sequencing, and path groups.

\textbf{Activity constraints:}

$$g^{(1)}_f: \quad x_3[r] \neq x_3[r'],\quad \forall r, r' \in \mathcal{R},\ r \neq r' \text{ s.t. } y_7[r] = y_7[r'] \quad \text{(unique vessel per role within an activity)},$$

$$g^{(2)}_f: \quad x_1[a] - \bigl(x_1[y_6[a]] + y_1[y_6[a]]\bigr) \geq 0,\quad \forall a \in \mathcal{A} : y_6[a] \neq \emptyset \quad \text{(precedence)}.$$

\textbf{Sequencing constraints:}

$$g^{(3)}_f: \quad \sum_{r' \in \mathcal{R}} x_4[r,r'] \leq 1,\quad \forall r \in \mathcal{R} \quad \text{(at most one successor)},$$

$$g^{(4)}_f: \quad \sum_{r \in \mathcal{R}} x_4[r,r'] \leq 1,\quad \forall r' \in \mathcal{R} \quad \text{(at most one predecessor)},$$

$$g^{(5)}_f: \quad x_4[r,r] \leq 0,\quad \forall r \in \mathcal{R} \quad \text{(no self-loops)},$$

$$g^{(6)}_f: \quad x_4[r,r'] = 1\ \Rightarrow\ x_3[r] = x_3[r'],\quad \forall r,r' \in \mathcal{R} \quad \text{(same vessel for consecutive roles)},$$

$$g^{(7)}_f: \quad x_4[r,r'] \leq 0,\quad \forall r,r' \in \mathcal{R} : y_7[r] = y_7[r'] \quad \text{(no intra-activity sequencing)},$$

$$g^{(8)}_f: \quad x_1[y_7[r]] - x_1[y_7[r']] < 0,\quad \forall r,r' \in \mathcal{R} : x_4[r,r'] = 1 \quad \text{(temporal ordering)},$$

$$g^{(9)}_f: \quad x_1[y_7[r']] - \Bigl(x_1[y_7[r]] + y_1[y_7[r]] + \theta\!\left(x_3[r],\ \ell^{\text{end}}(r),\ \ell^{\text{start}}(r'),\ y_{13}[x_3[r]]^{\max}\right)\Bigr) \geq 0$$
$$\forall r,r' \in \mathcal{R} : x_4[r,r'] = 1,\ x_3[r] = x_3[r'] \quad \text{(sufficient travel time)}.$$

\textbf{Path constraints:}

$$g^{(10)}_f: \quad \sum_{r \in \mathcal{R}} \mathbf{1}[x_3[r] = v] \cdot \mathbf{1}\!\left[\sum_{r' \in \mathcal{R}} x_4[r',r] = 0\right] = 1,\quad \forall v \in \mathcal{V} \quad \text{(unique sequence start per vessel)},$$

$$g^{(11)}_f: \quad \sum_{r \in \mathcal{R}} \mathbf{1}[x_3[r] = v] \cdot \mathbf{1}\!\left[\sum_{r' \in \mathcal{R}} x_4[r,r'] = 0\right] = 1,\quad \forall v \in \mathcal{V} \quad \text{(unique sequence end per vessel)},$$

$$g^{(12)}_f: \quad \sum_{r \in \mathcal{R}}\sum_{r' \in \mathcal{R}} x_4[r,r'] \cdot \mathbf{1}[x_3[r] = v] = \max\!\left(0,\ \sum_{r \in \mathcal{R}} \mathbf{1}[x_3[r] = v] - 1\right),\quad \forall v \in \mathcal{V} \quad \text{(path cardinality)},$$

$$g^{(13)}_f: \quad \bigl(x_4[r,r'] = 1 \wedge x_3[r] = v\bigr)\ \Rightarrow\ x_6[r,r'] \in y_{13}[v],\quad \forall v \in \mathcal{V},\ \forall r,r' \in \mathcal{R} \quad \text{(feasible sailing speeds)}.$$

\newpage
\section{DAS-CMOP Results}

\renewcommand{\arraystretch}{0.84}
\footnotesize
\begin{longtable}{@{} l r r r r r @{}}
\caption{
  Win frequencies for DAS-CMOP under Configuration C. 
}
\label{tab:dascmop_wins_c} \\
%
\toprule
\textbf{Problem} &
\textbf{IMAP-GA-II} &
\textbf{C-TAEA} &
\textbf{MOEA/D-CDP} &
\textbf{NSGA-II} &
\textbf{NSGA-III} \\
\midrule
\endfirsthead
%
\multicolumn{6}{c}{\tablename~\thetable{} \emph{(continued)}} \\[2pt]
\toprule
\textbf{Problem} &
\textbf{IMAP-GA-II} &
\textbf{C-TAEA} &
\textbf{MOEA/D-CDP} &
\textbf{NSGA-II} &
\textbf{NSGA-III} \\
\midrule
\endhead
%
\midrule
\multicolumn{6}{r}{\emph{Continued on next page}} \\
\endfoot
%
\bottomrule
\endlastfoot
%
\begin{filecontents*}{data/das_cmop_conflict_win_frequencies_latex.csv}
Problem,IMAP_GA,C_TAEA,MOEA_D_CDP,NSGA_II,NSGA_III,_group_break
DAS-CMOP1(1),\textbf{30},0,0,0,-,0
DAS-CMOP1(2),\textbf{30},0,0,0,-,0
DAS-CMOP1(3),\textbf{30},0,0,0,-,0
DAS-CMOP1(4),\textbf{30},0,0,0,-,0
DAS-CMOP1(5),\textbf{30},0,0,0,-,0
DAS-CMOP1(6),\textbf{30},0,0,0,-,0
DAS-CMOP1(7),\textbf{30},0,0,0,-,0
DAS-CMOP1(8),\textbf{30},0,0,0,-,0
DAS-CMOP1(9),\textbf{30},0,0,0,-,0
DAS-CMOP1(10),\textbf{30},0,0,0,-,0
DAS-CMOP1(11),\textbf{30},0,0,0,-,0
DAS-CMOP1(12),0,0,0,0,-,0
DAS-CMOP1(13),4,0,10,\textbf{16},-,0
DAS-CMOP1(14),12,0,4,\textbf{14},-,0
DAS-CMOP1(15),0,0,0,0,-,0
DAS-CMOP1(16),0,0,0,0,-,0
DAS-CMOP2(1),\textbf{30},0,0,0,-,1
DAS-CMOP2(2),\textbf{30},0,0,0,-,0
DAS-CMOP2(3),\textbf{30},0,0,0,-,0
DAS-CMOP2(4),\textbf{30},0,0,0,-,0
DAS-CMOP2(5),\textbf{30},0,0,0,-,0
DAS-CMOP2(6),\textbf{30},0,0,0,-,0
DAS-CMOP2(7),\textbf{30},0,0,0,-,0
DAS-CMOP2(8),\textbf{30},0,0,0,-,0
DAS-CMOP2(9),\textbf{30},0,0,0,-,0
DAS-CMOP2(10),\textbf{30},0,0,0,-,0
DAS-CMOP2(11),\textbf{30},0,0,0,-,0
DAS-CMOP2(12),\textbf{30},0,0,0,-,0
DAS-CMOP2(13),10,1,\textbf{12},7,-,0
DAS-CMOP2(14),9,0,8,\textbf{13},-,0
DAS-CMOP2(15),12,0,9,\textbf{13},-,0
DAS-CMOP2(16),\textbf{18},0,4,16,-,0
DAS-CMOP3(1),\textbf{30},0,0,0,-,1
DAS-CMOP3(2),\textbf{30},0,0,0,-,0
DAS-CMOP3(3),\textbf{30},0,0,0,-,0
DAS-CMOP3(4),\textbf{30},0,0,0,-,0
DAS-CMOP3(5),\textbf{30},0,0,0,-,0
DAS-CMOP3(6),\textbf{30},0,0,0,-,0
DAS-CMOP3(7),\textbf{30},0,0,0,-,0
DAS-CMOP3(8),\textbf{30},0,0,0,-,0
DAS-CMOP3(9),\textbf{30},0,0,0,-,0
DAS-CMOP3(10),\textbf{30},0,0,0,-,0
DAS-CMOP3(11),\textbf{30},0,0,0,-,0
DAS-CMOP3(12),\textbf{3},1,1,0,-,0
DAS-CMOP3(13),\textbf{14},0,5,11,-,0
DAS-CMOP3(14),5,1,11,\textbf{13},-,0
DAS-CMOP3(15),\textbf{21},0,3,14,-,0
DAS-CMOP3(16),\textbf{21},0,5,13,-,0
DAS-CMOP4(1),\textbf{30},0,0,0,-,1
DAS-CMOP4(2),\textbf{30},0,0,0,-,0
DAS-CMOP4(3),\textbf{30},0,0,0,-,0
DAS-CMOP4(4),\textbf{30},0,0,0,-,0
DAS-CMOP4(5),\textbf{30},0,0,0,-,0
DAS-CMOP4(6),\textbf{30},0,0,0,-,0
DAS-CMOP4(7),\textbf{30},0,0,0,-,0
DAS-CMOP4(8),\textbf{30},0,0,0,-,0
DAS-CMOP4(9),\textbf{30},0,0,0,-,0
DAS-CMOP4(10),\textbf{30},0,0,0,-,0
DAS-CMOP4(11),\textbf{30},0,0,0,-,0
DAS-CMOP4(12),6,14,3,\textbf{16},-,0
DAS-CMOP4(13),2,4,7,\textbf{17},-,0
DAS-CMOP4(14),8,3,2,\textbf{17},-,0
DAS-CMOP4(15),1,\textbf{19},10,0,-,0
DAS-CMOP4(16),1,10,\textbf{18},1,-,0
DAS-CMOP5(1),\textbf{30},0,0,0,-,1
DAS-CMOP5(2),\textbf{30},0,0,0,-,0
DAS-CMOP5(3),\textbf{30},0,0,0,-,0
DAS-CMOP5(4),\textbf{30},0,0,0,-,0
DAS-CMOP5(5),\textbf{30},0,0,0,-,0
DAS-CMOP5(6),10,0,\textbf{19},1,-,0
DAS-CMOP5(7),\textbf{30},0,0,0,-,0
DAS-CMOP5(8),\textbf{15},6,7,5,-,0
DAS-CMOP5(9),\textbf{30},0,0,0,-,0
DAS-CMOP5(10),8,0,\textbf{22},0,-,0
DAS-CMOP5(11),\textbf{30},0,0,0,-,0
DAS-CMOP5(12),3,12,2,\textbf{19},-,0
DAS-CMOP5(13),\textbf{21},0,9,0,-,0
DAS-CMOP5(14),\textbf{18},0,12,0,-,0
DAS-CMOP5(15),\textbf{13},2,8,7,-,0
DAS-CMOP5(16),9,7,5,\textbf{12},-,0
DAS-CMOP6(1),\textbf{30},0,0,0,-,1
DAS-CMOP6(2),\textbf{30},0,0,0,-,0
DAS-CMOP6(3),\textbf{30},0,0,0,-,0
DAS-CMOP6(4),\textbf{30},0,0,0,-,0
DAS-CMOP6(5),\textbf{30},0,0,0,-,0
DAS-CMOP6(6),\textbf{30},0,0,0,-,0
DAS-CMOP6(7),\textbf{30},0,0,0,-,0
DAS-CMOP6(8),\textbf{25},3,1,1,-,0
DAS-CMOP6(9),\textbf{30},0,0,0,-,0
DAS-CMOP6(10),14,0,0,\textbf{16},-,0
DAS-CMOP6(11),\textbf{30},0,0,0,-,0
DAS-CMOP6(12),\textbf{13},12,1,7,-,0
DAS-CMOP6(13),\textbf{26},0,0,4,-,0
DAS-CMOP6(14),\textbf{21},0,0,9,-,0
DAS-CMOP6(15),\textbf{16},1,2,12,-,0
DAS-CMOP6(16),6,3,4,\textbf{17},-,0
DAS-CMOP7(1),\textbf{30},0,0,-,0,1
DAS-CMOP7(2),\textbf{30},0,0,-,0,0
DAS-CMOP7(3),0,2,7,-,\textbf{21},0
DAS-CMOP7(4),\textbf{30},0,0,-,0,0
DAS-CMOP7(5),\textbf{30},0,0,-,0,0
DAS-CMOP7(6),\textbf{30},0,0,-,0,0
DAS-CMOP7(7),0,1,8,-,\textbf{21},0
DAS-CMOP7(8),\textbf{30},0,0,-,0,0
DAS-CMOP7(9),\textbf{30},0,0,-,0,0
DAS-CMOP7(10),\textbf{30},0,0,-,0,0
DAS-CMOP7(11),0,2,7,-,\textbf{21},0
DAS-CMOP7(12),\textbf{30},0,0,-,0,0
DAS-CMOP7(13),\textbf{30},0,0,-,0,0
DAS-CMOP7(14),\textbf{30},0,0,-,0,0
DAS-CMOP7(15),\textbf{30},0,0,-,0,0
DAS-CMOP7(16),\textbf{30},0,0,-,0,0
DAS-CMOP8(1),\textbf{30},0,0,-,0,1
DAS-CMOP8(2),\textbf{29},0,6,-,0,0
DAS-CMOP8(3),0,\textbf{21},3,-,6,0
DAS-CMOP8(4),\textbf{30},0,3,-,0,0
DAS-CMOP8(5),\textbf{20},0,2,-,8,0
DAS-CMOP8(6),\textbf{29},0,3,-,0,0
DAS-CMOP8(7),0,\textbf{21},2,-,7,0
DAS-CMOP8(8),\textbf{15},0,\textbf{15},-,0,0
DAS-CMOP8(9),\textbf{19},2,8,-,1,0
DAS-CMOP8(10),\textbf{30},0,4,-,0,0
DAS-CMOP8(11),0,\textbf{25},0,-,5,0
DAS-CMOP8(12),\textbf{18},0,11,-,1,0
DAS-CMOP8(13),\textbf{29},0,1,-,1,0
DAS-CMOP8(14),\textbf{19},0,11,-,0,0
DAS-CMOP8(15),\textbf{29},0,3,-,1,0
DAS-CMOP8(16),\textbf{20},0,10,-,0,0
DAS-CMOP9(1),\textbf{30},0,0,-,0,1
DAS-CMOP9(2),\textbf{19},0,13,-,0,0
DAS-CMOP9(3),2,1,\textbf{25},-,2,0
DAS-CMOP9(4),\textbf{22},0,10,-,0,0
DAS-CMOP9(5),\textbf{14},0,\textbf{14},-,2,0
DAS-CMOP9(6),\textbf{20},0,10,-,0,0
DAS-CMOP9(7),4,0,\textbf{21},-,5,0
DAS-CMOP9(8),\textbf{11},0,8,-,\textbf{11},0
DAS-CMOP9(9),\textbf{17},0,11,-,2,0
DAS-CMOP9(10),\textbf{23},0,7,-,0,0
DAS-CMOP9(11),1,3,\textbf{13},-,\textbf{13},0
DAS-CMOP9(12),9,0,\textbf{11},-,10,0
DAS-CMOP9(13),\textbf{18},0,\textbf{18},-,0,0
DAS-CMOP9(14),\textbf{19},0,6,-,5,0
DAS-CMOP9(15),\textbf{20},0,17,-,0,0
DAS-CMOP9(16),\textbf{19},0,4,-,7,0
\end{filecontents*}
\csvreader[
  head,
  late after line = \\,
]{data/das_cmop_conflict_win_frequencies_latex.csv}{}
{
  \ifnum\csvcolvii=1 \midrule \fi
  \csvcoli & \csvcolii & \csvcoliii & \csvcoliv & \csvcolv & \csvcolvi
}

\end{longtable}

\renewcommand{\arraystretch}{0.84}
\footnotesize
\begin{longtable}{@{} l r r r r r @{}}
\caption{
  Win frequencies for DAS-CMOP under Configuration A.
}
\label{tab:dascmop_wins_a} \\
%
\toprule
\textbf{Problem} &
\textbf{IMAP-GA-II} &
\textbf{C-TAEA} &
\textbf{MOEA/D-CDP} &
\textbf{NSGA-II} &
\textbf{NSGA-III} \\
\midrule
\endfirsthead
%
\multicolumn{6}{c}{\tablename~\thetable{} \emph{(continued)}} \\[2pt]
\toprule
\textbf{Problem} &
\textbf{IMAP-GA-II} &
\textbf{C-TAEA} &
\textbf{MOEA/D-CDP} &
\textbf{NSGA-II} &
\textbf{NSGA-III} \\
\midrule
\endhead
%
\midrule
\multicolumn{6}{r}{\emph{Continued on next page}} \\
\endfoot
%
\bottomrule
\endlastfoot
%
\begin{filecontents*}{data/das_cmop_no_conflict_win_frequencies_latex.csv}
Problem,IMAP_GA,C_TAEA,MOEA_D_CDP,NSGA_II,NSGA_III,_group_break
DAS-CMOP1(1),9,6,5,\textbf{10},-,0
DAS-CMOP1(2),6,7,2,\textbf{15},-,0
DAS-CMOP1(3),5,\textbf{12},7,6,-,0
DAS-CMOP1(4),6,11,1,\textbf{12},-,0
DAS-CMOP1(5),7,8,5,\textbf{10},-,0
DAS-CMOP1(6),6,7,5,\textbf{12},-,0
DAS-CMOP1(7),1,\textbf{23},0,6,-,0
DAS-CMOP1(8),3,\textbf{19},3,5,-,0
DAS-CMOP1(9),6,4,\textbf{10},\textbf{10},-,0
DAS-CMOP1(10),5,\textbf{10},7,8,-,0
DAS-CMOP1(11),1,\textbf{25},0,4,-,0
DAS-CMOP1(12),0,0,0,0,-,0
DAS-CMOP1(13),\textbf{14},0,10,6,-,0
DAS-CMOP1(14),\textbf{10},0,\textbf{10},\textbf{10},-,0
DAS-CMOP1(15),\textbf{3},2,0,1,-,0
DAS-CMOP1(16),0,0,0,0,-,0
DAS-CMOP2(1),\textbf{11},3,7,9,-,1
DAS-CMOP2(2),\textbf{19},0,3,8,-,0
DAS-CMOP2(3),\textbf{13},4,6,7,-,0
DAS-CMOP2(4),\textbf{17},2,8,3,-,0
DAS-CMOP2(5),\textbf{12},1,5,\textbf{12},-,0
DAS-CMOP2(6),9,1,9,\textbf{11},-,0
DAS-CMOP2(7),11,5,2,\textbf{12},-,0
DAS-CMOP2(8),\textbf{17},1,1,11,-,0
DAS-CMOP2(9),\textbf{21},0,5,4,-,0
DAS-CMOP2(10),9,1,8,\textbf{12},-,0
DAS-CMOP2(11),1,\textbf{21},0,8,-,0
DAS-CMOP2(12),\textbf{16},8,0,12,-,0
DAS-CMOP2(13),\textbf{12},0,7,11,-,0
DAS-CMOP2(14),8,0,7,\textbf{15},-,0
DAS-CMOP2(15),\textbf{15},0,2,\textbf{15},-,0
DAS-CMOP2(16),\textbf{18},0,1,17,-,0
DAS-CMOP3(1),\textbf{17},4,4,5,-,1
DAS-CMOP3(2),5,3,2,\textbf{20},-,0
DAS-CMOP3(3),\textbf{25},3,1,1,-,0
DAS-CMOP3(4),8,\textbf{13},2,7,-,0
DAS-CMOP3(5),\textbf{16},3,6,5,-,0
DAS-CMOP3(6),3,1,6,\textbf{20},-,0
DAS-CMOP3(7),\textbf{20},7,1,2,-,0
DAS-CMOP3(8),\textbf{16},1,3,10,-,0
DAS-CMOP3(9),\textbf{14},3,6,7,-,0
DAS-CMOP3(10),5,1,6,\textbf{18},-,0
DAS-CMOP3(11),\textbf{19},10,0,1,-,0
DAS-CMOP3(12),\textbf{14},3,2,7,-,0
DAS-CMOP3(13),\textbf{10},2,9,9,-,0
DAS-CMOP3(14),\textbf{11},0,\textbf{11},8,-,0
DAS-CMOP3(15),\textbf{14},4,5,\textbf{14},-,0
DAS-CMOP3(16),14,0,3,\textbf{20},-,0
DAS-CMOP4(1),\textbf{15},1,9,5,-,1
DAS-CMOP4(2),13,0,\textbf{19},0,-,0
DAS-CMOP4(3),0,9,\textbf{20},1,-,0
DAS-CMOP4(4),4,0,\textbf{28},0,-,0
DAS-CMOP4(5),\textbf{19},2,9,0,-,0
DAS-CMOP4(6),14,0,\textbf{17},0,-,0
DAS-CMOP4(7),0,8,\textbf{17},5,-,0
DAS-CMOP4(8),0,1,\textbf{26},3,-,0
DAS-CMOP4(9),10,0,\textbf{18},2,-,0
DAS-CMOP4(10),16,0,\textbf{17},0,-,0
DAS-CMOP4(11),0,\textbf{24},6,0,-,0
DAS-CMOP4(12),\textbf{17},6,6,5,-,0
DAS-CMOP4(13),\textbf{21},1,9,0,-,0
DAS-CMOP4(14),\textbf{21},1,14,2,-,0
DAS-CMOP4(15),5,0,4,\textbf{25},-,0
DAS-CMOP4(16),4,0,11,\textbf{15},-,0
DAS-CMOP5(1),0,0,0,\textbf{30},-,1
DAS-CMOP5(2),1,0,\textbf{16},13,-,0
DAS-CMOP5(3),0,0,0,\textbf{30},-,0
DAS-CMOP5(4),1,1,4,\textbf{24},-,0
DAS-CMOP5(5),0,\textbf{30},0,0,-,0
DAS-CMOP5(6),0,0,12,\textbf{18},-,0
DAS-CMOP5(7),1,3,4,\textbf{22},-,0
DAS-CMOP5(8),5,2,7,\textbf{16},-,0
DAS-CMOP5(9),0,0,0,\textbf{30},-,0
DAS-CMOP5(10),1,1,\textbf{15},13,-,0
DAS-CMOP5(11),8,\textbf{16},6,0,-,0
DAS-CMOP5(12),\textbf{13},6,2,\textbf{13},-,0
DAS-CMOP5(13),6,1,\textbf{14},9,-,0
DAS-CMOP5(14),8,1,\textbf{15},6,-,0
DAS-CMOP5(15),\textbf{16},9,4,3,-,0
DAS-CMOP5(16),\textbf{21},3,2,5,-,0
DAS-CMOP6(1),\textbf{15},1,14,0,-,1
DAS-CMOP6(2),\textbf{18},0,12,0,-,0
DAS-CMOP6(3),4,3,\textbf{22},1,-,0
DAS-CMOP6(4),8,0,\textbf{22},0,-,0
DAS-CMOP6(5),8,0,\textbf{19},3,-,0
DAS-CMOP6(6),\textbf{19},0,10,1,-,0
DAS-CMOP6(7),0,\textbf{14},\textbf{14},2,-,0
DAS-CMOP6(8),4,3,\textbf{13},10,-,0
DAS-CMOP6(9),\textbf{13},0,\textbf{13},4,-,0
DAS-CMOP6(10),\textbf{21},0,9,0,-,0
DAS-CMOP6(11),0,\textbf{28},2,0,-,0
DAS-CMOP6(12),15,\textbf{16},1,6,-,0
DAS-CMOP6(13),\textbf{27},0,3,0,-,0
DAS-CMOP6(14),\textbf{17},0,11,2,-,0
DAS-CMOP6(15),\textbf{24},0,5,1,-,0
DAS-CMOP6(16),\textbf{18},3,6,7,-,0
DAS-CMOP7(1),9,\textbf{10},4,-,7,1
DAS-CMOP7(2),\textbf{14},1,12,-,3,0
DAS-CMOP7(3),\textbf{12},10,2,-,6,0
DAS-CMOP7(4),\textbf{21},0,7,-,2,0
DAS-CMOP7(5),4,4,4,-,\textbf{18},0
DAS-CMOP7(6),\textbf{15},0,8,-,7,0
DAS-CMOP7(7),\textbf{8},7,\textbf{8},-,7,0
DAS-CMOP7(8),\textbf{12},0,7,-,11,0
DAS-CMOP7(9),3,10,5,-,\textbf{12},0
DAS-CMOP7(10),\textbf{19},1,8,-,2,0
DAS-CMOP7(11),\textbf{10},8,3,-,9,0
DAS-CMOP7(12),0,4,11,-,\textbf{15},0
DAS-CMOP7(13),13,0,\textbf{16},-,1,0
DAS-CMOP7(14),9,0,9,-,\textbf{12},0
DAS-CMOP7(15),\textbf{18},0,10,-,2,0
DAS-CMOP7(16),10,0,8,-,\textbf{12},0
DAS-CMOP8(1),\textbf{23},3,1,-,3,1
DAS-CMOP8(2),\textbf{26},0,3,-,1,0
DAS-CMOP8(3),2,\textbf{15},3,-,10,0
DAS-CMOP8(4),\textbf{22},0,7,-,1,0
DAS-CMOP8(5),\textbf{21},1,3,-,5,0
DAS-CMOP8(6),\textbf{23},0,5,-,3,0
DAS-CMOP8(7),6,\textbf{12},5,-,7,0
DAS-CMOP8(8),\textbf{16},0,9,-,5,0
DAS-CMOP8(9),\textbf{14},1,8,-,7,0
DAS-CMOP8(10),\textbf{19},0,5,-,6,0
DAS-CMOP8(11),1,\textbf{16},5,-,8,0
DAS-CMOP8(12),\textbf{15},0,12,-,3,0
DAS-CMOP8(13),\textbf{17},0,8,-,5,0
DAS-CMOP8(14),9,0,\textbf{15},-,6,0
DAS-CMOP8(15),\textbf{18},0,6,-,7,0
DAS-CMOP8(16),\textbf{12},0,9,-,9,0
DAS-CMOP9(1),\textbf{30},0,0,-,0,1
DAS-CMOP9(2),\textbf{23},0,6,-,1,0
DAS-CMOP9(3),4,\textbf{18},4,-,4,0
DAS-CMOP9(4),\textbf{20},2,8,-,1,0
DAS-CMOP9(5),5,\textbf{20},3,-,2,0
DAS-CMOP9(6),\textbf{22},1,8,-,0,0
DAS-CMOP9(7),\textbf{17},5,8,-,0,0
DAS-CMOP9(8),3,8,\textbf{12},-,7,0
DAS-CMOP9(9),3,\textbf{22},1,-,4,0
DAS-CMOP9(10),\textbf{24},1,4,-,1,0
DAS-CMOP9(11),\textbf{21},1,6,-,2,0
DAS-CMOP9(12),0,\textbf{17},8,-,5,0
DAS-CMOP9(13),13,5,\textbf{15},-,1,0
DAS-CMOP9(14),8,\textbf{10},7,-,5,0
DAS-CMOP9(15),13,5,\textbf{14},-,0,0
DAS-CMOP9(16),\textbf{11},10,5,-,4,0
\end{filecontents*}
\csvreader[
  head,
  late after line = \\,
]{data/das_cmop_no_conflict_win_frequencies_latex.csv}{}
{
  \ifnum\csvcolvii=1 \midrule \fi
  \csvcoli & \csvcolii & \csvcoliii & \csvcoliv & \csvcolv & \csvcolvi
}

\end{longtable}
\renewcommand{\arraystretch}{1.0}

\newpage
\section{DAS-CMaOP Results}
\label{appendix:dascmaop}

Tables~\ref{tab:dascmaop_c} and \ref{tab:dascmaop_a} report the complete per-problem win rates of the many-objective DAS-CMaOP campaign summarised in Section~\ref{subsec:results-dascmaop}: the mean win rate per DAS-CMaOP problem and objective count $m$, where each (problem, $m$) cell aggregates the twelve equality-free difficulty settings $\times$ ten seeds (120 runs), under Configurations C and A respectively. The best algorithm per row is in bold.

\renewcommand{\arraystretch}{0.84}
\footnotesize
\begin{longtable}{@{} l r r r r r @{}}
\caption{Mean win rate (\%) per DAS-CMaOP problem and objective count $m$ under Configuration C (conflicting preferences).}
\label{tab:dascmaop_c} \\
\toprule
\textbf{Problem} & $\mathbf{m}$ &
\textbf{IMAP-GA-II} &
\textbf{NSGA-III} &
\textbf{C-TAEA} &
\textbf{MOEA/D-CDP} \\
\midrule
\endfirsthead
\multicolumn{6}{c}{\tablename~\thetable{} \emph{(continued)}} \\[2pt]
\toprule
\textbf{Problem} & $\mathbf{m}$ &
\textbf{IMAP-GA-II} &
\textbf{NSGA-III} &
\textbf{C-TAEA} &
\textbf{MOEA/D-CDP} \\
\midrule
\endhead
\bottomrule
\endlastfoot
\begin{filecontents*}{data/dascmaop_conflict_mean_winrate_latex.csv}
Problem,m,IMAP_GA,NSGA_III,C_TAEA,MOEA_D_CDP,_group_break
DAS-CMaOP1,5,37.5,11.7,\textbf{38.3},2.5,0
DAS-CMaOP2,5,\textbf{73.3},5.8,10.0,5.0,0
DAS-CMaOP3,5,\textbf{56.7},5.8,27.5,1.7,0
DAS-CMaOP4,5,\textbf{71.7},4.2,10.0,9.2,0
DAS-CMaOP5,5,\textbf{70.0},9.2,10.0,3.3,0
DAS-CMaOP6,5,\textbf{77.5},2.5,10.0,1.7,0
DAS-CMaOP7,5,\textbf{60.0},4.2,16.7,10.8,0
DAS-CMaOP8,5,\textbf{80.8},1.7,4.2,5.8,0
DAS-CMaOP9,5,\textbf{75.8},3.3,10.8,1.7,0
DAS-CMaOP1,8,\textbf{44.2},11.7,23.3,0.0,1
DAS-CMaOP2,8,\textbf{89.2},0.0,1.7,0.0,0
DAS-CMaOP3,8,\textbf{67.5},0.0,15.8,0.0,0
DAS-CMaOP4,8,\textbf{83.3},0.0,0.0,1.7,0
DAS-CMaOP5,8,\textbf{90.8},0.0,1.7,0.8,0
DAS-CMaOP6,8,\textbf{80.8},0.0,2.5,0.0,0
DAS-CMaOP7,8,\textbf{79.2},1.7,2.5,0.8,0
DAS-CMaOP8,8,\textbf{83.3},0.0,1.7,0.0,0
DAS-CMaOP9,8,\textbf{80.8},0.0,1.7,0.8,0
DAS-CMaOP1,10,\textbf{43.3},4.2,35.8,0.0,1
DAS-CMaOP2,10,\textbf{90.0},0.0,0.0,1.7,0
DAS-CMaOP3,10,\textbf{62.5},2.5,16.7,1.7,0
DAS-CMaOP4,10,\textbf{79.2},4.2,0.0,0.0,0
DAS-CMaOP5,10,\textbf{79.2},3.3,0.8,0.0,0
DAS-CMaOP6,10,\textbf{83.3},0.0,0.0,0.0,0
DAS-CMaOP7,10,\textbf{80.8},0.8,1.7,0.0,0
DAS-CMaOP8,10,\textbf{81.7},0.0,2.5,0.0,0
DAS-CMaOP9,10,\textbf{73.3},7.5,2.5,0.0,0
\end{filecontents*}
\csvreader[
head,
late after line = \\,
]{data/dascmaop_conflict_mean_winrate_latex.csv}{}
{%
\ifnum\csvcolvii=1 \midrule \fi
\csvcoli & \csvcolii & \csvcoliii & \csvcoliv & \csvcolv & \csvcolvi
}
\end{longtable}

\begin{longtable}{@{} l r r r r r @{}}
\caption{Mean win rate (\%) per DAS-CMaOP problem and objective count $m$ under Configuration A (aligned preferences).}
\label{tab:dascmaop_a} \\
\toprule
\textbf{Problem} & $\mathbf{m}$ &
\textbf{IMAP-GA-II} &
\textbf{NSGA-III} &
\textbf{C-TAEA} &
\textbf{MOEA/D-CDP} \\
\midrule
\endfirsthead
\multicolumn{6}{c}{\tablename~\thetable{} \emph{(continued)}} \\[2pt]
\toprule
\textbf{Problem} & $\mathbf{m}$ &
\textbf{IMAP-GA-II} &
\textbf{NSGA-III} &
\textbf{C-TAEA} &
\textbf{MOEA/D-CDP} \\
\midrule
\endhead
\bottomrule
\endlastfoot
\begin{filecontents*}{data/dascmaop_no_conflict_mean_winrate_latex.csv}
Problem,m,IMAP_GA,NSGA_III,C_TAEA,MOEA_D_CDP,_group_break
DAS-CMaOP1,5,9.2,\textbf{50.8},12.5,17.5,0
DAS-CMaOP2,5,\textbf{81.7},4.2,4.2,4.2,0
DAS-CMaOP3,5,\textbf{65.8},25.8,0.0,0.0,0
DAS-CMaOP4,5,\textbf{45.8},19.2,5.0,25.8,0
DAS-CMaOP5,5,20.0,12.5,22.5,\textbf{37.5},0
DAS-CMaOP6,5,\textbf{65.8},9.2,13.3,3.3,0
DAS-CMaOP7,5,\textbf{32.5},27.5,4.2,27.5,0
DAS-CMaOP8,5,\textbf{61.7},15.0,10.0,5.8,0
DAS-CMaOP9,5,\textbf{53.3},18.3,18.3,1.7,0
DAS-CMaOP1,8,0.0,\textbf{64.2},10.0,4.2,1
DAS-CMaOP2,8,\textbf{80.0},4.2,3.3,4.2,0
DAS-CMaOP3,8,\textbf{50.8},27.5,4.2,0.8,0
DAS-CMaOP4,8,\textbf{55.0},10.8,3.3,18.3,0
DAS-CMaOP5,8,5.8,35.8,10.8,\textbf{39.2},0
DAS-CMaOP6,8,\textbf{41.7},9.2,10.0,22.5,0
DAS-CMaOP7,8,\textbf{48.3},6.7,12.5,16.7,0
DAS-CMaOP8,8,\textbf{67.5},5.8,6.7,3.3,0
DAS-CMaOP9,8,14.2,\textbf{28.3},12.5,\textbf{28.3},0
DAS-CMaOP1,10,23.3,\textbf{59.2},7.5,0.8,1
DAS-CMaOP2,10,\textbf{82.5},5.0,4.2,0.0,0
DAS-CMaOP3,10,\textbf{48.3},27.5,7.5,0.0,0
DAS-CMaOP4,10,\textbf{51.7},11.7,0.0,21.7,0
DAS-CMaOP5,10,16.7,28.3,3.3,\textbf{35.0},0
DAS-CMaOP6,10,\textbf{37.5},10.0,9.2,26.7,0
DAS-CMaOP7,10,\textbf{45.0},7.5,15.0,15.8,0
DAS-CMaOP8,10,\textbf{62.5},10.8,2.5,8.3,0
DAS-CMaOP9,10,19.2,\textbf{35.0},5.0,24.2,0
\end{filecontents*}
\csvreader[
head,
late after line = \\,
]{data/dascmaop_no_conflict_mean_winrate_latex.csv}{}
{
\ifnum\csvcolvii=1 \midrule \fi
\csvcoli & \csvcolii & \csvcoliii & \csvcoliv & \csvcolv & \csvcolvi
}
\end{longtable}
\normalsize
\renewcommand{\arraystretch}{1.0}

\newpage
\section{MO-VRPTW Results}

\renewcommand{\arraystretch}{0.84}
  \footnotesize
  \begin{longtable}{@{} l r r r r r @{}}
  \caption{%
    Win frequencies for MO-VRPTW under Configuration C (conflicting preferences).
    Values show the number of seeds (out of 3) in which each algorithm achieved
    the highest IMAP preference score on each Solomon instance.
    \textbf{IMAP-BRKGA} reports wins for the best-performing decoder variant per instance.
    Bold entries indicate the row winner; ties are both bolded.
  }
  \label{tab:vrptw_wins_c} \\
  \toprule
  \textbf{Instance} &
  \textbf{IMAP-BRKGA} &
  \textbf{BRKGA-WS} &
  \textbf{NSGA-II} &
  \textbf{NSGA-III} &
  \textbf{SPEA2} \\
  \midrule
  \endfirsthead
  \multicolumn{6}{c}{\tablename~\thetable{} \emph{(continued)}} \\[2pt]
  \toprule
  \textbf{Instance} &
  \textbf{IMAP-BRKGA} &
  \textbf{BRKGA-WS} &
  \textbf{NSGA-II} &
  \textbf{NSGA-III} &
  \textbf{SPEA2} \\
  \midrule
  \endhead
  \midrule
  \multicolumn{6}{r}{\emph{Continued on next page}} \\
  \endfoot
  \bottomrule
  \endlastfoot
  \begin{filecontents*}{data/vrptw_conflict_win_frequencies_latex.csv}
Instance,IMAP_BRKGA,BRKGA_WS,NSGA_II,NSGA_III,SPEA2,_group_break
C101,\textbf{3},0,0,0,0,0
C102,\textbf{3},0,0,0,0,0
C103,\textbf{2},0,1,0,0,0
C104,\textbf{3},0,0,0,0,0
C105,\textbf{3},0,0,0,0,0
C106,\textbf{3},0,0,0,0,0
C107,\textbf{3},0,0,0,0,0
C108,\textbf{3},0,0,0,0,0
C109,\textbf{3},0,0,0,0,0
C201,\textbf{3},0,\textbf{3},\textbf{3},\textbf{3},1
C202,\textbf{3},0,0,0,0,0
C203,\textbf{3},0,0,0,0,0
C204,\textbf{3},0,0,0,0,0
C205,\textbf{2},0,1,0,0,0
C206,\textbf{3},0,0,0,0,0
C207,\textbf{3},0,0,0,0,0
C208,\textbf{3},0,0,0,0,0
R101,\textbf{3},0,0,0,0,1
R102,\textbf{3},0,0,0,0,0
R103,\textbf{3},0,0,0,0,0
R104,\textbf{3},0,0,0,0,0
R105,\textbf{2},0,0,0,1,0
R106,\textbf{3},0,0,0,0,0
R107,\textbf{3},0,0,0,0,0
R108,\textbf{3},0,0,0,0,0
R109,\textbf{3},0,0,0,0,0
R110,\textbf{3},0,0,0,0,0
R111,\textbf{3},0,0,0,0,0
R112,\textbf{3},0,0,0,0,0
R201,\textbf{3},0,0,0,0,1
R202,\textbf{3},0,0,0,0,0
R203,\textbf{3},0,0,0,0,0
R204,\textbf{3},0,0,0,0,0
R205,\textbf{3},0,0,0,0,0
R206,\textbf{3},0,0,0,0,0
R207,\textbf{3},0,0,0,0,0
R208,\textbf{3},0,0,0,0,0
R209,\textbf{3},0,0,0,0,0
R210,\textbf{3},0,0,0,0,0
R211,\textbf{3},0,0,0,0,0
RC101,\textbf{3},0,0,0,0,1
RC102,\textbf{3},0,0,0,0,0
RC103,\textbf{2},0,0,0,1,0
RC104,\textbf{3},0,0,0,0,0
RC105,\textbf{2},0,1,0,0,0
RC106,\textbf{2},1,0,0,0,0
RC107,\textbf{3},0,0,0,0,0
RC108,\textbf{3},0,0,0,0,0
RC201,\textbf{3},0,0,0,0,1
RC202,\textbf{3},0,0,0,0,0
RC203,\textbf{3},0,0,0,0,0
RC204,\textbf{3},0,0,0,0,0
RC205,\textbf{3},0,0,0,0,0
RC206,\textbf{3},0,0,0,0,0
RC207,\textbf{3},0,0,0,0,0
RC208,\textbf{3},0,0,0,0,0
\end{filecontents*}
\csvreader[
    head,
    late after line = \\,
  ]{data/vrptw_conflict_win_frequencies_latex.csv}{}
  {
    \ifnum\csvcolvii=1 \midrule \fi
    \csvcoli & \csvcolii & \csvcoliii & \csvcoliv & \csvcolv & \csvcolvi
  }
  \end{longtable}
\renewcommand{\arraystretch}{1.0}

\newpage
\section{Relationship Between IMAP Optimality and Pareto Efficiency}
\label{appendix:pareto-imap}

This appendix clarifies when Pareto-front approximation followed by IMAP selection can recover the same decision point as direct IMAP-guided search, and when it cannot. The purpose is not to argue that Pareto-based methods are generally inappropriate. Rather, it identifies the subset of problems where the two views are compatible and the subset where direct preference-guided search is needed.

Throughout this appendix, \emph{IMAP-optimal} refers to the maximiser of the affine aggregate $\mathbf{A}(\cdot)$ of \eqref{eq:a-fine-aggregator} in the CMODP \eqref{eq:CMODP}, i.e.\ the solution $\mathbf{x}^\ast$ of \eqref{eq:imap_star}. Because the z-normalisation statistics $\mu_{k,i}$ and $\sigma_{k,i}$ are computed over an evaluation pool, IMAP scores, and hence IMAP optimality, are defined relative to a pool of alternatives; this is not a defect but a structural property of PFM-consistent aggregation, in which preference is contextual and rankings are relative to the set of alternatives under consideration \cite{wolfert2026uniqueaggregation}. The construction in \cite{wolfert2026uniqueaggregation} is stated for a finite set of alternatives; the framework in this paper operationalises it for arbitrary feasible sets by always scoring over a finite evaluation pool (the population during search, and the pooled representatives at evaluation time, Section~\ref{subsec:metrics}). Proposition~1 below is insensitive to this choice: its conclusion holds for every pool containing the compared solutions.

\paragraph{Proposition 1: aligned monotone preferences imply Pareto efficiency.}
Assume all objectives are stated in minimisation form. Suppose the preference model is \emph{aligned and monotone}: for every stakeholder $k$ and objective $i$, the preference function $p_{k,i}$ is non-increasing in the objective value, so that $f_i(\mathbf{x}) \leq f_i(\mathbf{y})$ implies $p_{k,i}(f_i(\mathbf{x})) \geq p_{k,i}(f_i(\mathbf{y}))$; and suppose additionally that every objective $i$ carries at least one preference function that is strictly decreasing in $f_i$ and has strictly positive weight $w_{k,i} > 0$. Then every IMAP-optimal feasible solution is Pareto-efficient with respect to the raw minimisation objectives, for any evaluation pool containing the solutions being compared.

\begin{proof}
    Assume, for contradiction, that an IMAP-optimal solution $\mathbf{x}^\ast$ is Pareto-dominated by another feasible solution $\mathbf{y}$ in the pool. Then $f_i(\mathbf{y}) \leq f_i(\mathbf{x}^\ast)$ for every objective and $f_j(\mathbf{y}) < f_j(\mathbf{x}^\ast)$ for at least one objective $j$. By monotonicity, $p_{k,i}(f_i(\mathbf{y})) \geq p_{k,i}(f_i(\mathbf{x}^\ast))$ for every $(k,i)$, and by the strictness assumption there is a pair $(k',j)$ with $w_{k',j} > 0$ and $p_{k',j}(f_j(\mathbf{y})) > p_{k',j}(f_j(\mathbf{x}^\ast))$.

    Consider the z-normalised scores over the pool. On the strictly improved dimension $(k',j)$, the pool contains two distinct preference values, so $\sigma_{k',j} > 0$; the z-transformation on that dimension is affine with positive scale $1/\sigma_{k',j}$ and therefore preserves the strict inequality: $z_{k',j}(\mathbf{y}) > z_{k',j}(\mathbf{x}^\ast)$. On every other dimension, either $\sigma_{k,i} > 0$, in which case the weak inequality is likewise preserved, or $\sigma_{k,i} = 0$, in which case every pool member has the same preference score on that dimension and the regularised z-score (Section~\ref{sec:problem_formulation}) is exactly $0$ for both solutions. Hence $z_{k,i}(\mathbf{y}) \geq z_{k,i}(\mathbf{x}^\ast)$ for all $(k,i)$, with strict inequality on a dimension of strictly positive weight, and since all weights are non-negative,
    \[
        \mathbf{A}(\mathbf{y}) = \sum_{k,i} w_{k,i}\, z_{k,i}(\mathbf{y})
        \;>\; \sum_{k,i} w_{k,i}\, z_{k,i}(\mathbf{x}^\ast) = \mathbf{A}(\mathbf{x}^\ast).
    \]
    This contradicts the IMAP optimality of $\mathbf{x}^\ast$. Therefore $\mathbf{x}^\ast$ is Pareto-efficient. Because the argument uses only per-dimension order preservation under positive scaling, it does not depend on which pool supplies $\mu_{k,i}$ and $\sigma_{k,i}$: the conclusion holds for every evaluation pool containing $\mathbf{x}^\ast$ and $\mathbf{y}$.
\end{proof}

This proposition explains the aligned-preference control used in the experiments: when preferences are aligned and monotone, the IMAP optimum lies on the Pareto front, so a sufficiently complete front approximation contains it among its candidates. One qualification is required before concluding that post-hoc IMAP selection over such a front recovers the same decision point as direct IMAP-guided search. Proposition~1 guarantees membership of the front, not invariance of the selection: because IMAP scores are pool-relative, restricting the evaluation pool from the feasible set to a front approximation changes the normalisation statistics, and the ranking among mutually non-dominated points may change with them. Post-hoc selection therefore returns the IMAP optimum \emph{relative to the front pool}. The two views agree exactly when both selections are scored in the same decision context, which is precisely what the common evaluation protocol of Section~\ref{subsec:metrics} enforces: all algorithms' representatives are pooled and scored together, so every method is judged by the same normalisation context. Under that protocol, the remaining empirical difference between the two approaches depends mainly on how well the Pareto method samples the relevant region of the front.

\paragraph{Proposition 2: the converse does not hold.}
Not every Pareto-efficient solution is IMAP-optimal.

\begin{proof}
    It is sufficient to give a counterexample. Consider a two-objective minimisation problem whose feasible set consists of three solutions $\mathbf{a}$, $\mathbf{b}$, and $\mathbf{c}$:
    \[
        F(\mathbf{a}) = (0,10), \qquad
        F(\mathbf{b}) = (10,0), \qquad
        F(\mathbf{c}) = (4,4).
    \]
    None of these solutions Pareto-dominates another. Solution $\mathbf{a}$ is best on the first objective but worse on the second, solution $\mathbf{b}$ is best on the second objective but worse on the first, and solution $\mathbf{c}$ is intermediate on both objectives. Hence all three solutions are Pareto-efficient.

    Now suppose there is one stakeholder, or an aligned group of stakeholders, with equal weights $w_1 = w_2 = 0.5$ and monotone decreasing preference functions
    \[
        p_1(f_1) = 100 - f_1, \qquad
        p_2(f_2) = 100 - f_2.
    \]
    Then the preference vectors are
    \[
        P(F(\mathbf{a})) = (100,90), \qquad
        P(F(\mathbf{b})) = (90,100), \qquad
        P(F(\mathbf{c})) = (96,96).
    \]
    Z-normalising each dimension over the pool $\{\mathbf{a}, \mathbf{b}, \mathbf{c}\}$ (each column has mean $95.\overline{3}$ and standard deviation $\approx 4.110$) gives, by symmetry,
    \[
        z(\mathbf{a}) = (1.136,\, -1.298), \qquad
        z(\mathbf{b}) = (-1.298,\, 1.136), \qquad
        z(\mathbf{c}) = (0.162,\, 0.162),
    \]
    so the aggregates are $\mathbf{A}(\mathbf{a}) = \mathbf{A}(\mathbf{b}) = -0.081$ and $\mathbf{A}(\mathbf{c}) = 0.162$. Solution $\mathbf{c}$, with the most balanced high preference scores, is the unique IMAP optimum. Thus $\mathbf{a}$ and $\mathbf{b}$ are Pareto-efficient but not IMAP-optimal.
\end{proof}

Pareto efficiency only states that no feasible solution improves all raw objectives simultaneously; it does not identify which efficient trade-off is most preferred by the stakeholder model. Therefore, even in the aligned monotone case, Pareto search still requires a post-hoc decision rule to choose one point from the front. This is the formal counterpart of the observation in \cite{wolfert2026uniqueaggregation} that Pareto-based group decision approaches can identify non-dominated fronts but do not by themselves produce a single, uniquely determined best-fit alternative.

\paragraph{Example: non-monotone or reversed preferences.}
Under non-monotone or reversed preferences, even the implication from raw-objective improvement to higher preference fails, so Pareto dominance with respect to the raw objectives no longer encodes the decision criterion. A minimal illustration suffices. Consider a single minimisation objective with two feasible solutions, $f(\mathbf{a}) = 0$ and $f(\mathbf{b}) = 10$, so that $\mathbf{a}$ is strictly better under the raw objective. Suppose the decision-maker's preference is target-seeking and peaks at $10$, for example
\[
    p(f) = 100 - (f-10)^2,
\]
giving $p(f(\mathbf{a})) = 0$ and $p(f(\mathbf{b})) = 100$: the IMAP-preferred solution is $\mathbf{b}$, the solution that is worse under the raw objective. Note that the quadratic acts on the \emph{physical performance value} before it enters the preference scale; it is an elicitation shape, not an operation on preference scores, and is therefore fully PFM-consistent — the affine-invariance requirements of \cite{Barzilai2010,wolfert2026uniqueaggregation} constrain operations on preference scores, not the form of the mapping from performance to preference. The same issue arises in multi-objective settings with reversed stakeholder preferences: a movement that improves a raw minimisation objective may increase one stakeholder's preference while decreasing another's. Raw-objective Pareto dominance is therefore not sufficient to represent the group-preference ordering when preferences are conflicting, reversed, or non-monotone; this is the regime targeted by Configuration~C in the experiments.

Hence Pareto-front approximation and IMAP-guided search agree only under additional assumptions: aligned monotone preferences (Proposition~1), adequate coverage of the relevant front region, and a shared evaluation pool for the final selection. Under conflicting, reversed, or non-monotone preference functions, the raw-objective Pareto relation no longer encodes the actual decision criterion, so direct preference-guided search is generally required.

\newpage

\section{Ablation of the Fixed Algorithmic Constants}
\label{appendix:ablation}

This appendix reports the one-at-a-time ablation of the four algorithmic constants held fixed throughout the experiments (Section~\ref{subsec:cv-pref}): the truncation threshold $\tau$, the feasibility weight $w_\mathrm{cv}$, and the tail parameters $k$ and $p$ of the constraint-violation preference function \eqref{eq:cv-pref}. Each constant is varied around the base configuration ($\tau = 40$, $w_\mathrm{cv} = 0.50$, $k = 2.0$, $p = 1.5$) with the other three held at base: $\tau \in \{0, 20, 60\}$ (where $\tau = 0$ effectively disables the discard), $w_\mathrm{cv} \in \{0.25, 0.75\}$, $k \in \{0.5, 8.0\}$, and $p \in \{1.0, 2.0\}$, giving nine variants plus the base.

Every variant is run with IMAP-GA-II under Configuration~C on DAS-CMOP1, DAS-CMOP3, and DAS-CMOP9 at difficulty levels 4, 8, and 14, with ten seeds per cell. To keep the 900-run study tractable, the ablation uses a reduced budget relative to the main campaign: population 40 and 1000 generations (40{,}000 evaluations, versus population 300 and 300{,}000 evaluations in Section~\ref{subsec:parameters}). Since all ten variants share this budget, the comparison among them is unaffected; absolute feasibility and win figures should not be compared with the main campaign tables. Evaluation follows the pooled protocol of Section~\ref{subsec:metrics}, with the ten variants taking the role of the compared algorithms: per (problem, difficulty, seed) cell, each variant's best feasible representative is pooled and scored by post-hoc IMAP aggregation over the substantive preference functions, and the highest-scoring variant wins the cell. Table~\ref{tab:ablation_constants} reports, per variant, the mean win share over the 90 cells, the mean feasibility rate of the final population, and the share of runs that returned no feasible representative at all.

\begin{table}[h]
    \centering
    \footnotesize
    \caption{One-at-a-time ablation of the fixed algorithmic constants (90 runs per variant: 3 problems $\times$ 3 difficulties $\times$ 10 seeds, Configuration C). Win shares are pool-relative across the ten variants, so the uniform reference level is 10\%; with 90 cells the two-standard-error band around it is approximately $[3.7, 16.3]$.}
    \label{tab:ablation_constants}
    \begin{tabular}{@{} l r r r @{}}
        \toprule
        \textbf{Variant} & \textbf{Win share (\%)} & \textbf{Feasibility (\%)} & \textbf{No feas.\ repr.\ (\%)} \\
        \midrule
        \begin{filecontents*}{data/ablation_summary_latex.csv}
variant,win,feas,norep,sep,
Base ($\tau{=}40$; $w_\mathrm{cv}{=}0.50$; $k{=}2.0$; $p{=}1.5$),11.1,88.3,0.0,0,
$\tau = 0$,12.2,81.6,0.0,1,
$\tau = 20$,16.7,91.9,0.0,0,
$\tau = 60$,11.1,89.6,0.0,0,
$w_\mathrm{cv} = 0.25$,3.3,43.2,56.7,1,
$w_\mathrm{cv} = 0.75$,3.3,98.1,0.0,0,
$k = 0.5$,12.2,88.7,0.0,1,
$k = 8.0$,11.1,89.2,0.0,0,
$p = 1.0$,12.2,91.2,0.0,1,
$p = 2.0$,6.7,79.9,11.1,0,
\end{filecontents*}
\csvreader[
        head,
        late after line = \\,
        ]{data/ablation_summary_latex.csv}{}
        {\csvcoli & \csvcolii & \csvcoliii & \csvcoliv}
        \bottomrule
    \end{tabular}
\end{table}

Two findings follow. First, no variant outperforms the base configuration beyond sampling noise: all win shares except those of the $w_\mathrm{cv}$ extremes and $p = 2.0$ lie within the two-standard-error band around the uniform level, and the single value at its upper edge ($\tau = 20$) is consistent with chance under nine simultaneous comparisons. The framework is thus insensitive to $\tau$ and $k$ over the ranges tested, including the complete removal of the discard step ($\tau = 0$), which costs a few points of feasibility but does not change the outcome distribution significantly.

Second, the only statistically significant deviations are degradations at the extremes, and they fail in exactly the directions the design rationale of Section~\ref{subsec:cv-pref} predicts. Halving the feasibility weight ($w_\mathrm{cv} = 0.25$) understeers: mean feasibility collapses to 43.2\% and 56.7\% of runs end without any feasible representative. Raising it ($w_\mathrm{cv} = 0.75$) oversteers: feasibility rises to 98.1\%, but the win share falls to 3.3\% because search effort is spent on feasibility that could have been spent on preference improvement. Steepening the penalty tail ($p = 2.0$) reproduces the premature-convergence failure mode anticipated in Section~\ref{subsec:cv-pref}: feasibility drops on the hardest difficulty and 11.1\% of runs return no feasible representative. The base configuration therefore sits between two empirically confirmed failure modes rather than at an arbitrary point, and none of the neighbouring configurations offers a significant improvement, so the constants are retained unchanged for all main experiments. The scope of this study is deliberately limited to a DAS-CMOP subset with one-at-a-time variation; interaction effects between the constants and their behaviour on the combinatorial domains remain future work.

\end{document}